\documentclass[10pt]{IEEEtran}

\usepackage{graphics} 
\usepackage{epsfig} 
\usepackage{amsmath} 
\usepackage{amssymb}  
\usepackage{epstopdf}
\usepackage[amsmath,thmmarks]{ntheorem}
\usepackage[compress]{cite}
\usepackage{fix-cm}
\usepackage{enumerate}
\usepackage{bbm}
\usepackage{setspace}
\usepackage{dsfont}
\usepackage{color}
\usepackage{physics}
\usepackage[normalem]{ulem}
\usepackage [english]{babel}
\usepackage [autostyle, english = american]{csquotes}

\MakeOuterQuote{"}
\ifCLASSOPTIONcompsoc
\usepackage[caption=false,font=normalsize,labelfont=sf,textfont=sf]{subfig}
\else
\usepackage[caption=false,font=footnotesize]{subfig}
\fi
\theoremheaderfont{\normalfont\bfseries}
\theorembodyfont{\normalfont}
\theoremseparator{:}

\theoremnumbering{arabic}
\newtheorem{Theorem}{Theorem}

\theoremnumbering{arabic}
\newtheorem{Proposition}{Proposition}

\theoremnumbering{arabic}
\newtheorem{Corollary}{Corollary}

\theoremnumbering{arabic}
\theoremsymbol{\ensuremath{\square}}

\theoremsymbol{\ensuremath{\square}}
\theoremnumbering{arabic}

\theoremsymbol{\ensuremath{\square}}
\theoremnumbering{arabic}

\theoremsymbol{\ensuremath{\square}}
\theoremseparator{:}
\theoremnumbering{arabic}

\theoremsymbol{\ensuremath{\square}}
\theoremnumbering{arabic}
\newtheorem{Definition}{Definition}

\theoremstyle{break}

\theoremheaderfont{\itshape}
\theorembodyfont{\normalfont}
\theoremstyle{nonumberplain}
\theoremseparator{:}
\theoremsymbol{\ensuremath{\square}}
\newtheorem{Proof}{Proof}

\begin{document}
\title{Linear-Quadratic Time-Inconsistent Mean-Field Type Stackelberg Differential Games: Time-Consistent Open-Loop Solutions} 

\author{
Jun Moon,~ \IEEEmembership{Member, IEEE} and Hyun Jong Yang,~ \IEEEmembership{Member, IEEE}
\thanks{This research was supported in part by the National Research Foundation of Korea (NRF) Grant funded by the Ministry of Science and ICT, Korea (NRF-2017R1E1A1A03070936, NRF-2017R1A5A1015311), and in part by Institute for Information \& communications Technology Promotion (IITP) grant funded by the Korea government (MSIT), Korea (No. 2018-0-00958).}
\thanks{Jun Moon is with the School of Electrical and Computer Engineering, University of Seoul, Seoul 02504, South Korea; email:  \texttt{jmoon12@uos.ac.kr}.}
\thanks{Hyun Jong Yang is with the School of Electrical and Computer Engineering, Ulsan National Institute of Science and Technology (UNIST), Ulsan 44919, South Korea; email: \texttt{hjyang@unist.ac.kr}.}
}
\maketitle

\begin{abstract}
In this paper, we consider the linear-quadratic time-inconsistent mean-field leader-follower Stackelberg stochastic differential game with an adapted open-loop information structure. Given a controlled linear stochastic differential equation (SDE), the quadratic objective functionals of the leader and the follower include conditional expectations of state and control  (mean field) variables. In addition, the cost parameters could be general nonexponential discounting depending on the initial time. As stated in the existing literature, these two general settings of the objective functionals induce time inconsistency in the optimal solutions. Given an arbitrary control of the leader, we first obtain the follower's (time-consistent) equilibrium control and its state feedback representation in terms of the nonsymmetric coupled Riccati differential equations (RDEs) and the backward SDE. This provides the rational behavior of the follower, characterized by the forward-backward SDE (FBSDE). We then obtain the leader's explicit (time-consistent) equilibrium control and its state feedback representation in terms of the nonsymmetric coupled RDEs, where the constraint of the leader's problem is the FBSDE induced by the follower's rational behavior. With the solvability of the nonsymmetric coupled RDEs, the equilibrium controls of the leader and the follower constitute the time-consistent Stackelberg equilibrium of the paper. Finally, the numerical examples are provided to check the solvability of the nonsymmetric coupled RDEs of the leader and the follower.
\end{abstract}
\begin{IEEEkeywords}
time-inconsistent stochastic control problem, equilibrium control, Stackelberg differential games, mean-field stochastic systems.
\end{IEEEkeywords}

\IEEEpeerreviewmaketitle

\section{Introduction}\label{Section_1}

One of the most important features of (deterministic and stochastic) optimal control that is typically not emphasized  is time consistency, which states that the optimal solution obtained with respect to the initial time $t \in [0,T]$ and the initial state $x(t) = x$ (equivalently the initial pair $(t,x)$) remains optimal when it is restricted on $[\tau, T]$ with the initial condition $x(\tau)$ (equivalently the initial pair $(\tau,x(\tau))$ with $\tau \in (t,T]$ \cite{Hu_SICON_2012, Yong_AMAS_2012, Yong_AMS_2017, Yan_Inbook_2019}.\footnote{As mentioned in \cite{Hu_SICON_2012}, it can also be stated that the optimal control viewed at present will also be optimal when it is viewed at a later time.} The time-consistent property stems from Bellman's principle of optimality (dynamic programming), which formulates the optimal control problem as a family of the initial pair $(t,x)$. In fact, the dynamic programming approach provides a theoretical foundation for optimal control theory (and differential games) by characterizing the value function in terms of the initial pair $(t,x)$, which leads to the Hamilton-Jacobi-Bellman (HJB) equation \cite{Liberzon_Book, Yong_book, Touzi_Book}.

We expect that \emph{all} (deterministic and stochastic) optimal control problems hold the time-consistent property. However, in many situations, the time-consistency fails. Specifically, as discussed in the examples of \cite[Examples 1.1 and 1.2]{Yong_AMS_2017}, \cite[Section 2]{Yong_MCRF_2012} and \cite{Yong_MCRF_2011, Yan_Inbook_2019}, this is due to the fact that
\begin{enumerate}
\item[(1)] the cost parameters in the objective functional could be general nonexponential discounting parameters that are dependent on the initial time;
\item[(2)] the conditional expectations of state and control variables could be included in the state equation (e.g. stochastic differential equation, SDE) and the objective functional.
\end{enumerate}
Under the above two general settings, \emph{time-inconsistent (deterministic and stochastic) optimal control problems} can be formulated. Regarding (1), in \cite{Strotz_RES_1956, Loewenstein_QJE_1992, Krusell_Econometrica_2003, Karp_JPE_2005, Ekeland_arXiv_2006, Ekeland_MFE_2008, Bjork_2010, Ekeland_MFE_2010, Ekeland_SIFM_2012, Marin_EJOR_2010}, time-inconsistent (deterministic and stochastic) optimal control problems with general (nonexponential) discounting parameters including quasi-exponential discounting, hyperbolic discounting, and quasi-geometric discounting were considered. The purpose of introducing general nonexponential discounting parameters in the objective functional is to capture the various discounted preferences of the optimizer (or user) with respect to the initial time and state. As for (2), the time-consistency fails since the law of iterated expectations cannot be applied to obtain a dynamic programming equation \cite{Bjork_2010}. In fact, the conditional expectations of state and control variables are known as \emph{mean field} variables. The motivation to include mean-field variables in the state equation (e.g. SDE) and the objective functional is to analyze macroscopic behavior of large-scale interacting particle systems in engineering, biology and economics, and to consider mean-variance portfolio optimization in various mathematical finance applications \cite{Hu_SICON_2012, 4303232, Weintraub_2008, Bensoussan_book_2013, Benssousan_AMO_2016, Bjork_2010, Yong_MCRF_2012, Yong_AMS_2017, Wei_SICON_2017, Yong_SICON_2013_MF, Hu_SICON_2017, Lasry, Moon_TAC, Moon_Automatica_2018, Ni_TAC_2018, TAC_2019_accepted_MK, Elie_AOR_2019}.

There are various definitions for quantifying \enquote{optimality} in (deterministic and stochastic) time-inconsistent optimal control problems. First, one could search for a \emph{precommitted} optimal solution, which is optimal only for a prescribed initial pair $(t,x)$. This definition corresponds to the standard optimal solution, which can be characterized by a usual variational approach. Note that the precommitted optimal solution is time-inconsistent; hence, it is not implementable in the sense that at a later time $\tau > t$, we have to search for a second (precommitted) optimal solution \cite{Bjork_2010, Ekeland_SIFM_2012}. Second, instead of a precommitted optimal solution, an \emph{equilibrium control} can be sought, which is local with respect to spike variation but time consistent \cite[Definition 2.1]{Hu_SICON_2017}.\footnote{It is called the equilibrium control, since it is an equilibrium with respect to all his admissible controls in the future \cite{Hu_SICON_2012}.} The equilibrium control concept is related to an (adapted for the stochastic case) open-loop control, which can be obtained by using the (stochastic) maximum principle. One advantage of this definition is that it is mathematically rigorous, and in some situations the equilibrium control can be expressed as a feedback representation using the decoupling method \cite{Hu_SICON_2012, Yong_AMS_2017, Hu_SICON_2017}. The results of \cite{Hu_SICON_2012, Hu_SICON_2017} were extended to the jump diffusion case in \cite{Sun_JOTA_2019, Alia_MCRF_2019}. The last definition corresponds to the closed-loop or state-feedback type control. This can be characterized by solving the HJB equation obtained by discretizing the corresponding optimal control problem, which is closely related to the multi-person differential game \cite{Ekeland_arXiv_2006, Ekeland_MFE_2008, Ekeland_MFE_2010, Ekeland_SIFM_2012, Yong_MCRF_2011, Yong_MCRF_2012, Yong_AMAS_2012, Yong_AMS_2017, Bjork_2010, Wei_SICON_2017, Bjork_FS_2017, Wei_MCRF_2017}. In addition, the mixed equilibrium solution concept was used in \cite{Ni_SICON_2019}, and the Markovian framework for time-inconsistent linear-quadratic problems was developed in \cite{Wang_AMO_2018}.

We note that the aforementioned references correspond to a class of time-inconsistent optimal control problems, where there is one single control (decision maker) in the state equation and the objective functional. Within this formulation, hierarchical decision-making analysis between players cannot be considered. Then it is natural to extend the earlier results for time-inconsistent optimal control problems to the leader-follower game framework, subject to a prescribed decision-making hierarchy. It should be noted that the time-inconsistent leader-follower game has not been considered in the existing literature, and this problem is addressed in our paper (see the problem formulation and the summary of the main results of the paper in Section \ref{Section_1_1}).

The class of leader-follower differential games is also known as \emph{Stackelberg} differential games \cite{Basar2, Yong_SICON_2002, Bensoussan_SICON_2015}. The leader holds a dominating position; the leader chooses and then announces his optimal strategy by considering the rational behavior of the follower. Under this hierarchical setting, 
the leader's optimal solution and the follower's rational behavior constitute a Stackelberg equilibrium. Classical (deterministic and stochastic) Stackelberg differential games have been studied extensively in the literature; see \cite{Basar2, 1101986, 1100307, 1101999, Freiling, Yong_SICON_2002, Bensoussan_SICON_2015, Moon_Automatica_2018} and the references therein. Note that depending on the open-loop or closed-loop information structure between the leader and the follower, the Stackelberg game has to be treated differently.\footnote{See \cite{Basar2} and \cite{Bensoussan_SICON_2015} for definitions and discussions of open-loop and closed-loop information structures in Stackelberg differential games, and the approaches with respect to different information structures.} In particular, for the adapted open-loop case, the (stochastic) Stackelberg game can be analyzed by applying the stochastic maximum principle to the following: (i) the follower's problem given an arbitrary strategy of the leader and (ii) the leader's problem, where the constraint is the follower's rational behavior characterized by the forward-backward SDE (FBSDE) from (i) \cite{Basar2, Yong_SICON_2002, Bensoussan_SICON_2015}. There are wide ranges of applications for Stackelberg games including engineering, economics and biology; see \cite{Chen_TAC_1972, Srikant_large, Chu_IEC_2014, 6464552, Moon_Automatica_2018, Lin_TAC_2018} and the references therein.

\subsection{Problem Statement and  Main Results of the Paper}\label{Section_1_1}
In this paper, we consider the \emph{linear-quadratic (LQ) time-inconsistent mean-field Stackelberg stochastic differential game} for the leader and the follower. The adapted open-loop information structure is adopted in the sense that the follower chooses his optimal decision after the leader announces his optimal strategy over the entire horizon \cite{Yong_SICON_2002, Bensoussan_SICON_2015}. In the problem setting, a linear stochastic differential equation (SDE) controlled by the leader and the follower is given, where their control variables are also included in the diffusion term of the SDE.\footnote{This implies that the magnitude of the stochastic noise is controlled by the leader and the follower, which can be viewed as \emph{multiplicative noise} of the system \cite{Ugrinovskii_SICON_1999, Moon_TAC_2019}.} The objective functionals of the leader and the follower are quadratic, where the conditional expectations of state and control variables (mean field) are nonlinearly included, and the cost parameters could be general nonexponential discounting, depending on the initial time. As mentioned, these two general settings of the objective functionals induce the time inconsistency of optimal solutions for the leader and the follower.

Our main results of the paper can be summarized as follows:
\begin{enumerate}[(i)]
	\item we obtain the follower's equilibrium control\footnote{\label{footnote_1}As mentioned, an equilibrium control is time consistent; see Definition \ref{Definition_2} or \cite[Definition 4.1]{Yong_AMS_2017} and \cite[Definition 2.1]{Hu_SICON_2017}.} that is a function of an arbitrary leader's control (see Proposition \ref{Proposition_1}). Then we obtain its state feedback representation in terms of the nonsymmetric coupled Riccati differential equations (RDEs) and the backward stochastic differential equation (BSDE) (see Theorem \ref{Theorem_1}). Note that the follower's equilibrium control induces the rational behavior of the follower that is characterized by the forward-backward SDE (FBSDE);
	\item the explicit leader's equilibrium control\textsuperscript{\ref{footnote_1}} is obtained, where the constraint of the leader's equilibrium control problem is the rational behavior of the follower that is the FBSDE from (i) (see Theorem \ref{Theorem_2}). We obtain the state feedback representation of the leader's equilibrium control in terms of the nonsymmetric coupled RDEs via the generalized decoupling technique (see Theorem \ref{Theorem_3});
	\item with the solvability of the nonsymmetric coupled RDEs in (i) and (ii), the results of (i) and (ii) constitute the time-consistent (adapted open-loop) Stackelberg equilibrium for the leader and the follower (see Corollary \ref{Corollary_1}); 
	\item numerical examples are provided to check the solvability of the nonsymmetric coupled RDEs in (i) and (ii).
\end{enumerate}

Note that the problem formulation and the results of the paper can be viewed as extensions of those in \cite{Yong_SICON_2002, Lin_TAC_2018} to the time-inconsistent problem, and those in \cite{Yong_SICON_2013_MF, Yong_AMS_2017} to the Stackelberg game framework. In \cite{Lin_TAC_2018}, the LQ mean-field Stackelberg game was considered, where the corresponding Stackelberg equilibrium is precommitted, i.e., it is time inconsistent. The extensions of \cite{Yong_SICON_2002, Lin_TAC_2018} to the time-inconsistent setting are not trivial, since the approach for characterizing the (time-consistent) equilibrium control is completely different from that of the the precommitted (time-inconsistent) optimal solution as discussed in \cite{Hu_SICON_2012, Yong_AMS_2017, Hu_SICON_2017}. On the other hand, in \cite{Yong_AMS_2017} the time-inconsistent mean-field control problem was studied, where the (adapted open-loop time-consistent) equilibrium control was obtained. Note that the extension of \cite{Yong_AMS_2017} to the Stackelberg game is also challenging, since in (ii) in the preceding list, we need to solve the leader's time-inconsistent stochastic optimal control problem with the FBSDE constraint induced by the follower. We mention that the time-inconsistent  stochastic optimal control problem with the FBSDE constraint has not been studied in the existing literature.

The paper is organized as follows. The problem formulation is stated in Section \ref{Section_2}. The equilibrium control problems of the follower and the leader are considered in Sections \ref{Section_3} and \ref{Section_4}, respectively. Numerical examples are presented in Section \ref{Section_5}. The concluding remarks are given in Section \ref{Section_6}. 

\textit{Notation}: Let $\mathbb{R}^n$ be the $n$-dimensional Euclidian space, and $\mathbb{S}^n$ the set of $n \times n$ dimensional symmetric matrices. For $x \in \mathbb{R}^n$, $x^\top$ denotes its transpose. Let $\langle x,y \rangle$ be the inner product and $|x| := \langle x, x \rangle^{1/2}$ for $x, y \in \mathbb{R}^n$. Let $|x|^2_{S} = x^\top S x$ for $ x \in \mathbb{R}^n$ and $S \in \mathbb{S}^n$. For $X \in \mathbb{S}^n$, let $X > 0$ (resp. $X \geq 0$) be a positive semidefinite (resp. positive definite) matrix. $I$ denotes an identity matrix with an appropriate dimension. Let $L^\infty([t,T],\mathbb{R}^n)$ be the set of $\mathbb{R}^n$-valued functions with $\|f\|_\infty = \sup_{s \in [t,T]} |f(s)| < \infty$ for $f \in L^\infty([t,T],\mathbb{R}^n)$. Let $(\Omega, \mathcal{F}, \mathbb{P},\{\mathcal{F}_t\}_{t \geq 0})$ be a complete filtered probability space on which a one dimensional standard Brownian motion $B(\cdot)$ is defined, where $\mathbb{F} = \{\mathcal{F}_t\}_{t \geq 0}$ is a natural filtration generated by the Brownian motion augmented by all the $\mathbb{P}$-null sets in $\mathcal{F}$. Let $\mathbb{E}$ be the mathematical expectation operator and $\mathbb{E}_t[\cdot] = \mathbb{E}[\cdot|\mathcal{F}_t]$. 
Let $\mathcal{S}_\mathbb{F}([t,T],\mathbb{R}^n)$ be the set of $\mathbb{R}^n$-valued $\mathbb{F}$-adapted stochastic processes such that for $x \in \mathcal{S}_\mathbb{F}([t,T],\mathbb{R}^n)$, $x$ is continuous and satisfies $\mathbb{E}_t[\sup_{s \in [t,T]} |x(s)|^2 ] < \infty$. For $p \geq 1$, let $\mathcal{L}_\mathbb{F}^p([t,T],\mathbb{R}^n)$ be the set of $\mathbb{R}^n$-valued $\mathbb{F}$-adapted stochastic processes such that for $x \in \mathcal{L}_\mathbb{F}^p([t,T],\mathbb{R}^n)$, $x$ satisfies $\mathbb{E}_t[\int_t^T |x(s)|^p \dd s] <\infty$.

\section{Problem Statement}\label{Section_2}

We consider the following stochastic differential equation (SDE) driven by Brownian motion:
\begin{align}
\label{eq_1}
\begin{cases}
\dd x(s) = \Bigl [ A(s)x(s) + B_1(s)u(s) + B_2(s) v(s) \Bigr ]\dd s \\
 + \Bigl [ C(s) x(s) + D_1(s) u(s) + D_2(s) v(s) \Bigr ] \dd B(s),~ s \in [t,T] \\
x(t) = x_0,
\end{cases}
\end{align}
where $x \in \mathbb{R}^n$ is state with the initial condition $x(t) = x_0$ and $t \in [0,T]$, $u \in \mathbb{R}^{m_1}$ is control of the leader and $v \in \mathbb{R}^{m_2} $ is control of the follower. In (\ref{eq_1}), $A,C : [0,T] \rightarrow \mathbb{R}^{n \times n}$ and $B_i, D_i : [0,T] \rightarrow \mathbb{R}^{n \times m_i}$, $i=1,2$, are deterministic coefficient matrices with $A,C \in L^\infty([0,T],\mathbb{R}^n)$ and $B_i,D_i \in L^\infty([0,T],\mathbb{R}^{m_i})$ for $i=1,2$. The sets of admissible controls of the leader is defined as follows:
\begin{align*}
& \mathcal{U}[t,T] \\
&= \{u:[t,T] \times \Omega \times \mathbb{R}^n \rightarrow \mathbb{R}^{m_1}~|~ u \in \mathcal{L}_{\mathbb{F}}^2([t,T],\mathbb{R}^{m_1})\}. 
\end{align*}
The set of admissible controls for the follower, $\mathcal{V}[t,T]$, is defined in a similar way. Then for any $u \in \mathcal{U}[t,T]$ and $v \in \mathcal{V}[t,T]$, there exists a unique solution of (\ref{eq_1}) that satisfies $\mathbb{E}[\sup_{s \in [t,T]} |x(s)|^2]< \infty$, i.e., $x \in \mathcal{S}_\mathbb{F}([t,T],\mathbb{R}^n)$ \cite[Chapter 1, Theorem 6.14]{Yong_book} (see also \cite[Theorem 2.2]{Touzi_Book}).\footnote{The assumption of the one-dimensional Brownian motion in (\ref{eq_1}) is only for notational convenience, and we can easily extend the results of the paper to the multi-dimensional Brownian motion case.} Let $x^{u,v}$ be the state process in (\ref{eq_1}) controlled by $u \in \mathcal{U}[t,T]$ and $v \in \mathcal{V}[t,T]$. 

The objective functional of the leader is given by
\begin{align}
\label{eq_2}
	& J_1(t,x_0;u, v)  \\
	& = \mathbb{E}_t \Biggl [ \int_t^T \bigl [ |x(s)|^2_{Q_1(s,t)} + |\mathbb{E}_t [x(s)]|^2_{\bar{Q}_1(s,t)} + |u(s)|^2_{R_1(s,t)} \nonumber \\
	&~~~  + |\mathbb{E}_t [u(s)]|^2_{\bar{R}_1 (s,t)} \bigr ] \dd s + |x(T)|^2_{M_1(t)} + |\mathbb{E}_t [x(T)]|^2_{\bar{M}_1(t)} \Biggr ], \nonumber 
\end{align}
and the objective functional of the follower is as follows
\begin{align}
\label{eq_3}
	& J_2(t,x_0;u, v)  \\
	& = \mathbb{E}_t \Biggl [ \int_t^T \bigl [ |x(s)|^2_{Q_2(s,t)} + |\mathbb{E}_t [x(s)]|^2_{\bar{Q}_2(s,t)} + |v(s)|^2_{R_2(s,t)} \nonumber \\
	&~~~  + |\mathbb{E}_t [v(s)]|^2_{\bar{R}_2(s,t)} \bigr ] \dd s + |x(T)|^2_{M_2(t)} + |\mathbb{E}_t [x(T)]|^2_{\bar{M}_2(t)} \Biggr ]. \nonumber
\end{align}
In (\ref{eq_2}) and (\ref{eq_3}), the cost parameters hold $Q_i, \bar{Q}_i \in L^\infty([0,T] \times [0,T],\mathbb{S}^n)$, $M_i, \bar{M}_i \in L^\infty([0,T], \mathbb{S}^n)$ and $R_i, \bar{R}_i \in L^\infty([0,T] \times [0,T],\mathbb{S}^{m_i})$ for $i=1,2$. It is assumed that
\begin{align}
\label{eq_4}
\begin{cases}
Q_i(s,t) \geq 0,~ Q_i(s,t) + \bar{Q}_i(s,t) \geq 0 \\
M_i(t) \geq 0,~ M_i(t) \geq 0, M_i(t) + \bar{M}_i(t) \geq 0 \\
R_i(s,t) > 0,~ R_i(s,t) + \bar{R}_i(s,t) > 0 \\
s \in [t,T],~ t \in [0,T],~i=1,2.
\end{cases}
\end{align}
In (\ref{eq_2}) and (\ref{eq_3}), the cost parameters depend on the initial time $t$. In addition, not only the state and control variables, but also their conditional expectations are included \emph{nonlinearly} in the objective functionals in (\ref{eq_2}) and (\ref{eq_3}). The conditional expectation terms in (\ref{eq_2}) and (\ref{eq_3}) are also known as \emph{mean field} of state and control variables \cite{Yong_SICON_2013_MF, Yong_AMS_2017, Ni_TAC_2018, Lin_TAC_2018, TAC_2019_accepted_MK}.

The interaction between the leader and the follower can be stated as follows. The leader chooses and announces his optimal solution to the follower by considering the rational reaction of the follower. The follower then determines his optimal solution by responding to the optimal solution of the leader. Under this setting, the problem can be solved in a reverse way \cite{Yong_SICON_2002, Bensoussan_SICON_2015}. Specifically,
\begin{enumerate}
\setlength{\itemindent}{0.5em}
\item[(S.1)] solve the follower's optimal control problem with arbitrary $x_0 \in \mathbb{R}^n$ and control of the leader $u \in \mathcal{U}[t,T]$, that is, minimize $J_2(t,x_0;u,v)$ over $v \in \mathcal{V}[t,T]$ subject to (\ref{eq_1}) for any $x_0 \in \mathbb{R}^n$ and $u \in \mathcal{U}[t,T]$, where its solution is denoted by $\bar{v} \in \mathcal{V}[t,T]$;
\item[(S.2)] obtain the leader's optimal solution, denoted by $\bar{u} \in \mathcal{U}[t,T]$, by minimizing $J_1(t,x_0;u,\bar{v})$ over $u \in \mathcal{U}[t,T]$ subject to (\ref{eq_1}) with $v$ replaced by the follower's optimal solution $\bar{v}$ obtained from (i).
\end{enumerate}
In general, the optimal solution of the follower from (S.1), $\bar{v} \in \mathcal{V}[t,T]$, depends on an arbitrary $u \in \mathcal{U}[t,T]$ and the initial condition $x_0 \in \mathbb{R}^n$. In fact, (S.1) above characterizes the rational reaction behavior of the follower, which is the optimization constraint of (S.2).

Due to the hierarchy between the leader and the follower mentioned above, the problem considered in the paper can be referred to as the (adapted open-loop) linear-quadratic (LQ) mean-field Stackelberg differential game \cite{Basar2, Yong_SICON_2002, Bensoussan_SICON_2015}. We also note that if the solutions to (S.1) and (S.2) exist, then the pair $(\bar{u}, \bar{v})$ constitutes the (adapted open-loop) Stackelberg equilibrium \cite{Yong_SICON_2002, Bensoussan_SICON_2015} (see also Definition \ref{Definition_1} below).

Now, suppose that the cost parameters in (\ref{eq_2}) and (\ref{eq_3}) do not depend on the initial time, i.e., for $i=1,2$,
\begin{align}
\label{eq_5}
\begin{cases}
	Q_i(s,t) = Q_i(s),~ \bar{Q}_i(s,t) = \bar{Q}_i(s) \\
	R_i(s,t) = R_i(s),~ \bar{R}_i(s,t) =\bar{R}_i(s) \\
	M_i(t) = M_i,~ \bar{M}_i(t) = \bar{M}_i.
\end{cases}
\end{align}
Then with (\ref{eq_4}) and (\ref{eq_5}), the (adapted open-loop) Stackelberg equilibrium was obtained in \cite[Theorems 3.1-3.3]{Lin_TAC_2018}, where its precise definition is given as follows \cite[Definition 2.1]{Lin_TAC_2018}:

\begin{Definition}\label{Definition_1}	
The pair $(\bar{u},\bar{v}) \in \mathcal{U}[t,T] \times \mathcal{V}[t,T]$ constitutes an (adapted open-loop) \emph{Stackelberg equilibrium}, and $\bar{x} := x^{\bar{u},\bar{v}} \in \mathcal{S}_{\mathbb{F}}([t,T],\mathbb{R}^n)$ is the corresponding state process if the following conditions hold.
\begin{enumerate}
\item[(i)] There exists a measurable map $\alpha : \mathcal{U}[t,T] \times \mathbb{R}^n \rightarrow \mathcal{V}[t,T]$ such that for any $u \in \mathcal{U}[t,T]$,
\begin{align*}	
& J_2(t,x_0;u,\alpha[u,x_0]) 
= \inf_{ v \in \mathcal{V}[t,T]}  J_2(t,x_0;u,v);
\end{align*}
\item[(ii)]There exists a control $\bar{u} \in \mathcal{U}[t,T]$ such that
\begin{align*}	
& J_1(t,x_0;\bar{u},\alpha[\bar{u},x_0])  
= \inf_{ u \in \mathcal{U}[t,T]}  J_1(t,x_0;u,\alpha[u,x_0]);
\end{align*}
\item[(iii)]$\bar{v} = \alpha[\bar{u}, x_0]$ with the state process $\bar{x} = x^{\bar{u},\bar{v}} = x^{\bar{u},\alpha[\bar{u}, x_0]}$.
\end{enumerate}
\end{Definition}
The results in \cite{Lin_TAC_2018} were obtained via the variational method. This approach can also be applied to the problem of the paper (that is without assuming (\ref{eq_5})). However, in this case the corresponding (adapted open-loop) Stackelberg equilibrium (see Definition \ref{Definition_1}) would be \emph{time inconsistent}\footnote{Specifically, under (\ref{eq_4}), with the objective functionals in (\ref{eq_2}) and (\ref{eq_3}), the Stackelberg equilibrium in Definition \ref{Definition_1} can be obtained by using the results in \cite[Theorems 3.1-3.3]{Lin_TAC_2018}. However, the corresponding Stackelberg equilibrium is time inconsistent due to (F.1) and (F.2) \cite{Yong_AMS_2017}.} in the sense that the Stackelberg solution at time $t$ may not be the Stackelberg solution at any $\mathbb{F}$-stopping time $\tau \in (t,T]$ \cite{Yong_MCRF_2012, Hu_SICON_2012, Yong_MCRF_2011, Yong_AMS_2017, Hu_SICON_2017, Ni_TAC_2018, Wei_SICON_2017, Lin_TAC_2018, Yan_Inbook_2019}. This is due to the fact that
\begin{enumerate}
\setlength{\itemindent}{0.5em}
\item[(F.1)] the cost parameters in (\ref{eq_2}) and (\ref{eq_3}) (which are $Q_i(s,t)$, $\bar{Q}_i(s,t)$, $R_i(s,t)$, $\bar{R}_i(s,t)$, $M_i(t)$ and $\bar{M}_i(t)$, $i=1,2$) could be general nonexponential discounting parameters that are dependent on the initial time $t$ \cite{Yong_MCRF_2011, Yong_AMS_2017, Ni_TAC_2018, Alia_MCRF_2019};
\item[(F.2)] the mean field (the conditional expectations) of state and control variables are included in the objective functionals in a \emph{nonlinear} way \cite{Hu_SICON_2012, Yong_MCRF_2011, Yong_AMS_2017, Hu_SICON_2017, Wei_SICON_2017}.
\end{enumerate}
In view of the above discussion and due to (F.1) and (F.2), the problem of the paper can be regarded as the \emph{linear-quadratic (LQ) time-inconsistent mean-field Stackelberg differential game} (\texttt{Problem LQ-TI-MF-SDG}). We also mention that the time-inconsistent Stackelberg equilibrium studied in \cite{Lin_TAC_2018} with Definition \ref{Definition_1} is closely related to the \emph{precommitted} optimal solutions of the leader and the follower, since they are optimal only when viewed at the initial time \cite{Bjork_2010, Hu_SICON_2012, Yong_AMS_2017}.

Now, given the time-inconsistent nature of \texttt{Problem LQ-TI-MF-SDG} as indicated in the preceding section, the objective of our paper is to obtain the \emph{time-consistent} Stackelberg solution of \texttt{Problem LQ-TI-MF-SDG}. Hence, it is necessary to modify the notion of \enquote{optimality} in Definition \ref{Definition_1} using the the time-consistent equilibrium solution concept used in the time-inconsistent stochastic optimal control problems studied in \cite{Yong_MCRF_2011, Hu_SICON_2012, Yong_AMS_2017, Wei_SICON_2017, Ni_TAC_2018, Hu_SICON_2017} (see \cite[Definition 4.1]{Yong_AMS_2017} and \cite[Definition 2.1]{Hu_SICON_2017}). Before stating its specific definition, we introduce the following \enquote{infinitesimally} perturbed control of the leader via spike variation: Given a control $u^* \in \mathcal{U}[t,T]$, for almost all $t \in [0,T)$, $\epsilon > 0$ and any $ u \in \mathcal{U}[t,T]$, define
\begin{align}
\label{eq_6}
u^\epsilon(s) = u^*(s)\mathbbm{1}_{s \in [t+\epsilon,T]} + u(s) \mathbbm{1}_{s \in [t,t+\epsilon)},~ s \in [t,T],
\end{align}
where $\mathbbm{1}$ is an indicator function. Similarly, given $\alpha[u,x_0]:\mathcal{U}[t,T] \times \mathbb{R}^n \rightarrow \mathcal{V}[t,T]$ with a fixed $u \in \mathcal{U}[t,T]$, for almost all $t \in [0,T)$, $\epsilon > 0$ and any $v \in \mathcal{V}[t,T]$ with $s \in [t,T]$, define
\begin{align}	
\label{eq_7}
\alpha^\epsilon[u,x_0](s)  &= 	\alpha[u,x_0](s)\mathbbm{1}_{s \in [t+\epsilon,T]}  + v(s) \mathbbm{1}_{s \in [t,t+\epsilon)}.
\end{align}


\begin{Definition}\label{Definition_2}
	The pair $(u^*,v^*) \in \mathcal{U}[t,T] \times \mathcal{V}[t,T]$ constitutes a \emph{time-consistent (adapted open-loop) Stackelberg equilibrium} of \texttt{Problem LQ-TI-MF-SDG} for the leader and the follower, and $x^* := x^{u^*,v^*} \in \mathcal{S}_\mathbb{F}([t,T],\mathbb{R}^n)$ with $x^*(t) = x_0$ is corresponding the equilibrium state process if for any initial condition $x_0$, the following conditions hold:
\begin{enumerate}[(i)]
	\item For any given $u \in \mathcal{U}[t,T]$ and almost all $t \in [0,T)$, there are measurable map $\alpha:\mathcal{U}[t,T] \times \mathbb{R}^n \rightarrow \mathcal{V}[t,T]$ and the corresponding state process $x^{u,\alpha[u,x_0]}$ such that	
	\begin{align*}
	& \liminf_{\epsilon \downarrow 0} \frac{1}{\epsilon} \Bigl ( J_2(t,x_0;u,\alpha^\epsilon[u,x_0]) \\
	& ~~~~~~ - J_2(t,x_0;u,\alpha[u,x_0])  \Bigr )\geq 0,
	\end{align*}
	where $\alpha^\epsilon$ is defined in (\ref{eq_7}). For the pair $(x^{u,\alpha[u,x_0]},\alpha)$, $\alpha$ is called the \emph{equilibrium control of the follower} under an arbitrary $u \in \mathcal{U}[t,T]$ of the leader, and $x^{u,\alpha[u,x_0]}$ is the corresponding \emph{equilibrium state process}.
	\item For almost all $t \in [0,T)$, there are $u^* \in \mathcal{U}[t,T]$ and the corresponding state process $x^{u^*,\alpha[u^*, x_0]}$ such that 
	\begin{align*}
	& \liminf_{\epsilon \downarrow 0} \frac{1}{\epsilon} \Bigl ( J_1(t,x_0;u^\epsilon,\alpha[u^\epsilon, x_0]) \\
	&~~~~~~ - J_1(t,x_0;u^*,\alpha[u^*, x_0]) \Bigr ) \geq 0,
	\end{align*}
	where $u^\epsilon$ is defined in (\ref{eq_6}) and $\alpha$ is the equilibrium control of the follower in (i). For the pair $(x^{u^*,\alpha[u^*, x_0]},u^*)$, $u^*$ is called the \emph{equilibrium control of the leader} and $x^{u^*,\alpha[u^*, x_0]}$ is the corresponding \emph{equilibrium state process}.
	\item[(iii)] $v^* = \alpha[u^*,x_0]$ with the equilibrium state process $x^* = x^{u^*,v^*} = x^{u^*,\alpha[u^*, x_0]}$.
\end{enumerate}
\end{Definition}

Note that Definition \ref{Definition_2} is local in an infinitesimal sense. In (i) and (ii) of Definition \ref{Definition_2}, the two equilibrium pairs, $(x^{u,\alpha[u,x_0]},\alpha)$ and $(x^{u^*,\alpha[u^*, x_0]},u^*)$, can be viewed as the one-player (time-consistent) equilibrium solution for the time-inconsistent stochastic optimal control problem \cite{Hu_SICON_2012, Yong_AMS_2017, Hu_SICON_2017}. Hence, Definition \ref{Definition_2}(i) implies that for any $u \in \mathcal{U}[t,T]$, the follower plays a game at any time $t$ against all his admissible controls in the future. 
The same argument applies to the leader's case with Definition \ref{Definition_2}(ii). 



\section{Follower's Equilibrium Control} \label{Section_3}
This section considers the follower's equilibrium control problem in the sense of Definition \ref{Definition_2}(i). 
The follower's problem is denoted by \texttt{Problem LQ-FEC} (LQ follower's equilibrium control problem). We state the following result:

\begin{Proposition}\label{Proposition_1}
Consider \texttt{Problem LQ-FEC} with (\ref{eq_4}). For any $u \in \mathcal{U}[t,T]$, suppose that the pair $(x^*,v^*)$, where $v^* := \alpha : \mathcal{U}[t,T] \times \mathbb{R}^n \rightarrow \mathcal{V}[t,T]$ and $x^* := x^{u,\alpha[u,x_0]} = x^{u,v^*}$ with $x^*(t) = x_0$, 
is the state-control pair of the follower. Consider the backward stochastic differential equation (BSDE):
\begin{align}
\label{eq_8}
\begin{cases}
\dd p(s,t) = - \bigl [ A^\top(s) p(s,t) + C^\top(s) q(s,t) + Q_2(s,t) x^*(s) \\
	~~  + \bar{Q}_2(s,t) \mathbb{E}_t[x^*(s)] \bigr ] \dd s + q(s,t) \dd B(s),~ s \in [t,T] \\
p(T,t) = M_2(t) x(T) + \bar{M}_2(t) \mathbb{E}_t[x(T)].
\end{cases}
\end{align}
Suppose that $(p,q) \in \mathcal{L}_{\mathbb{F}}^2([0,T] \times [0,T],\mathbb{R}^n) \times \mathcal{L}_{\mathbb{F}}^2([0,T] \times [0,T],\mathbb{R}^n)$ is the solution to the BSDE in (\ref{eq_8}), which is continuous in both $s \in [t,T]$ and $t \in [0,T]$. Then the corresponding unique (open-loop) equilibrium control of the follower satisfies for $t \in [0,T]$, with $\widehat{R}_2(s,t) = R_2(s,t) + \bar{R}_2(s,t)$,
\begin{align}
\label{eq_9}
  \widehat{R}_2(t,t) v^*(t) + B_2^\top(t) p(t,t) + D_2^\top(t) q(t,t) = 0.
\end{align}
In particular, for $t \in [0,T]$, we have
\begin{align}	
\label{eq_10}
v^*(t) = - \widehat{R}_2^{-1}(t,t) [B_2^\top(t) p(t,t) + D_2^\top(t) q(t,t)].
\end{align}
\begin{Proof}
The result follows from \cite[Proposition 4.2]{Yong_AMS_2017} (see also \cite[Theorem 3.2]{Hu_SICON_2012} or \cite[Theorem 3.5]{Hu_SICON_2017}) with the perturbed control defined in (\ref{eq_7}). Also, by using the approach in \cite[Theorem 3.5]{Hu_SICON_2017}, we can show that (\ref{eq_9}) is the necessary and sufficient condition for the equilibrium control of the follower, which implies the uniqueness. This completes the proof.
\end{Proof}
\end{Proposition}

In view of Proposition \ref{Proposition_1}, we have the forward-backward stochastic differential equation (FBSDE) with the optimality condition in (\ref{eq_10}) (note that $^*$ notation is dropped in $x$): 
\begin{align}
\label{eq_11}
\begin{cases}
	\dd x(s) = \bigl [ A(s) x(s) + B_1(s) u(s) \\
	~~~ + B_2(s) v^*(s) \bigr ] \dd s + \bigl [ C(s) x(s) + D_1(s) u(s) \\
	~~~ + D_2(s) v^*(s) \bigr ] \dd B(s),~ s \in [t,T] \\
	\dd p(s,t) = - \bigl [ A^\top(s) p(s,t) + C^\top(s) q(s,t)  \\
	~~~  + Q_2(s,t) x(s) + \bar{Q}_2(s,t) \mathbb{E}_t[x(s)] \bigr ] \dd s \\
	~~~ + q(s,t) \dd B(s),~ s \in [t,T]\\
	x(t) = x_0\\
	p(T,t) = M_2(t) x(T) + \bar{M}_2(t) \mathbb{E}_t[x(T)],
\end{cases}	
\end{align}
where $v^*$ satisfies the optimality condition in (\ref{eq_9}) that is given below for convenience
\begin{align*}
	\widehat{R}_2(t,t) v^*(t) + B_2^\top(t) p(t,t) + D_2^\top(t) q(t,t) = 0.
\end{align*}
This corresponds to the rational behavior of the follower for an arbitrary control of the leader.

We now obtain the state feedback representation of the open-loop equilibrium control of the follower in (\ref{eq_10}). By applying the Four-Step Scheme\footnote{This approach is quite similar to decoupling between adjoint and state equations used in LQ optimal control problems \cite{Yong_book, Liberzon_Book}.} used in \cite{Yong_book, Ma_book, Moon_TAC_2018, Yong_SICON_2002}, we consider the following transformation to decouple forward and backward parts of the SDE:
\begin{align}
\label{eq_12}
p(s,t) = P(s,t) x(s) + Z(s,t) \mathbb{E}_t[x(s)] + h(s),
\end{align}
where the explicit expressions of $P$ and $Z$ will be obtained later. Here, it is assumed that $P,Z : [0,T] \times [0,T] \rightarrow \mathbb{R}^{n \times n}$, and $h$ is the first component of the $n$-dimensional BSDE that will be characterized later. Note that in view of the terminal condition of (\ref{eq_8}), we must have $P(T,t) = M_2(t)$, $Z(T,t) = \bar{M}_2(t)$ and $h(T) = 0$. Let $\widehat{P}(s,t) = P(s,t) + Z(s,t)$. 

By using It\^o's formula, we have
\begin{align}
\label{eq_13}
& \dd p(s,t) \\
&= - \Bigl [ A^\top(s) p(s,t) + C^\top(s) q(s,t) + Q_2(s,t) x(s) \nonumber \\
	&~~~  + \bar{Q}_2(s,t) \mathbb{E}_t[x(s)] \Bigr ] \dd s + q(s,t) \dd B(s) \nonumber \\
	& = \frac{\dd P(s,t)}{\dd s} x(s) \dd s + 	\frac{\dd Z(s,t)}{\dd s} \mathbb{E}_t[x(s)] \dd s \nonumber \\
&~~~  + P(s,t)\Bigl [A(s) x(s) + B_1(s) u(s) + B_2(s) v^*(s) \Bigr ] \dd s \nonumber \\
&~~~ + P(s,t) \Bigl [ C(s) x(s) + D_1(s) u(s) + D_2(s) v^*(s) \Bigr ] \dd B(s) \nonumber \\
&~~~ + Z(s,t) \Bigl [ A(s) \mathbb{E}_t[x(s)] + B_1(s) \mathbb{E}_t[u(s)] \nonumber \\
&~~~ + B_2(s) \mathbb{E}_t[v^*(s)] \Bigr ] \dd s  + k(s) \dd s + l(s) \dd B(s), \nonumber
\end{align}
where $\dd h(s) = k(s) \dd s + l(s) \dd B(s)$. Then
\begin{align}	
\label{eq_14}
 q(s,t) & = P(s,t) \bigl [ C(s) x(s) \\
 &~~~ + D_1(s) u(s) + D_2(s) v^*(s) \bigr ] + l(s). \nonumber 
\end{align}

By substituting (\ref{eq_12}) and (\ref{eq_14}) into the optimality condition in (\ref{eq_9}), we have
\begin{align*}
	& \widehat{R}_2(t,t) v^*(t) + B_2^\top(t) \bigl [ \widehat{P}(t,t)x(t) + h(t) \bigr ]  + D_2^\top(t) P(t,t) \\
	& \times \bigl [ C(t) x(t) + D_1(t) u(t) + D_2(t) v^*(t) \bigr ] + D_2^\top(t) l(t) = 0,
\end{align*}
which implies
\begin{align}
\label{eq_15}
& v^*(t) = \alpha[u,x_0](t)  \\
& = - \bigl [ \widehat{R}_2(t,t) + D_2^\top(t) P(t,t)D_2(s) \bigr]^{-1} \nonumber \\
&~~~ \times  \bigl [ \bigl (B_2^\top(t) \widehat{P}(t,t) + D_2^\top(t) P(t,t) C(t) \bigr ) x(t) \nonumber \\
&~~~ +  \bigl ( B_2^\top(t) h(t) + D_2^\top(t) P(t,t) D_1(t) u(t) + D_2^\top(t) l(t) \bigr ) \bigr ]. \nonumber
\end{align}
$v^*$ in (\ref{eq_15}) is the feedback representation of the follower's equilibrium control in (\ref{eq_10}) that is obtained by the transformations in (\ref{eq_12}) and (\ref{eq_14}). Note that (\ref{eq_15}) depends on an arbitrary control of the leader as in Definition \ref{Definition_2}(i).

From (\ref{eq_14}) and (\ref{eq_15}), we have
\begin{align}
\label{eq_16}
& q(s,t) = P(s,t) C(s) x(s) + P(s,t) D_1(s) u(s) + l(s)  \\
&  - P(s,t) D_2(s) \bigl [ \widehat{R}_2(s,s) + D_2^\top(s) P(s,s)D_2(s) \bigr]^{-1} \nonumber \\
& \times \bigl [ B_2^\top(s) \widehat{P}(s,s) + D_2^\top(s) P(s,s) C(s) \bigr ] x(s) \nonumber \\
& - P(s,t) D_2(s) \bigl [ \widehat{R}_2(s,s) + D_2^\top(s) P(s,s)D_2(s) \bigr]^{-1} \nonumber \\
& \times \bigl [ B_2^\top(s) h(s) + D_2^\top(s) P(s,s) D_1(s) u(s) + D_2^\top(s) l(s) \bigr ]. \nonumber
\end{align}
By substituting (\ref{eq_12}), (\ref{eq_15}) and (\ref{eq_16}) into (\ref{eq_13}), we have
\begin{align*}
& \dd p(s,t) \\
&=  - \Bigl [ A^\top(s) P(s,t) x(s) +  A^\top(s)  Z(s,t) \mathbb{E}_t[x(s)] \\
&~~~ +  A^\top(s)  h(s)  + C^\top(s) P(s,t) C(s) x(s) \\
&~~~ + C^\top(s) P(s,t) D_1(s) u(s) + C^\top(s)  l(s) \\
&~~~ - C^\top(s) P(s,t) D_2(s) \\
&~~~ \times \bigl [ \widehat{R}_2(s,s) + D_2^\top(s) P(s,s)D_2(s) \bigr]^{-1} \\
&~~~ \times \bigl [ B_2^\top(s) \widehat{P}(s,s) + D_2^\top(s) P(s,s) C(s) \bigr ] x(s) \\
&~~~ - C^\top(s) P(s,t) D_2(s) \\
&~~~ \times \bigl [ \widehat{R}_2(s,s) + D_2^\top(s) P(s,s)D_2(s) \bigr]^{-1} \\
&~~~ \times \bigl [ B_2^\top(s) h(s) + D_2^\top(s) P(s,s) D_1(s) u(s) \\
&~~~ + D_2^\top(s) l(s) \bigr ]  + Q_2(s,t) x(s) \\
&~~~ + \bar{Q}_2(s,t) \mathbb{E}_t[x(s)] \Bigr ] \dd s + q(s,t) \dd B(s),
\end{align*}
and
\begin{align*}
& \dd p(s,t)  = \frac{\dd P(s,t)}{\dd s} x(s) \dd s + 	\frac{\dd Z(s,t)}{\dd s} \mathbb{E}_t[x(s)] \dd s \\
& ~~~ + P(s,t) A(s) x(s) \dd s + P(s,t) B_1(s) u(s) \dd s \\
&~~~ - P(s,t) B_2(s) \Bigl [ \widehat{R}_2(s,s) + D_2^\top(s) P(s,s)D_2(s) \Bigr]^{-1} \\
&~~~ \times \Bigl [ B_2^\top(s) \widehat{P}(s,s) + D_2^\top(s) P(s,s) C(s) \Bigr ] x(s) \dd s \\
&~~~ - P(s,t) B_2(s) \Bigl [ \widehat{R}_2(s,s) + D_2^\top(s) P(s,s)D_2(s) \Bigr]^{-1} \\
&~~~ \times \Bigl [ B_2^\top(s) h(s) + D_2^\top(s) P(s,s) D_1(s) u(s) \\
&~~~ + D_2^\top(s) l(s) \Bigr ]\dd s  + P(s,t) \\
&~~~ \times \Bigl [ C(s) x(s) + D_1(s) u(s) + D_2(s) v^*(s) \Bigr ] \dd B(s) \\
&~~~ + Z(s,t) A(s) \mathbb{E}_t[x(s)] \dd s + Z(s,t) B_1(s) \mathbb{E}_t[u(s)]  \dd s \\
&~~~ - Z(s,t) B_2(s) \Bigl [ \widehat{R}_2(s,s) + D_2^\top(s) P(s,s)D_2(s) \Bigr]^{-1} \\
&~~~ \times \Bigl [ B_2^\top(s) \widehat{P}(s,s) + D_2^\top(s) P(s,s) C(s) \Bigr ] \mathbb{E}_t[x(s)] \dd s \\
&~~~ - Z(s,t) B_2(s) \Bigl [ \widehat{R}_2(s,s) + D_2^\top(s) P(s,s)D_2(s) \Bigr]^{-1} \\
&~~~ \times \Bigl [ B_2^\top(s) \mathbb{E}_t[h(s)] + D_2^\top(s) P(s,s) D_1(s) \mathbb{E}_t[u(s)] \\
&~~~ + D_2^\top(s) \mathbb{E}_t[l(s)] \Bigr ]\dd s  + k(s,t) \dd s + l(s) \dd B(s).
\end{align*}

We compare the coefficients in the above two equalities. Then we can easily show that $P$ and $Z$ have to satisfy the following nonsymmetric coupled Riccati differential equations (RDEs):
\begin{align}
\label{eq_17}
\begin{cases}
- \frac{\dd P(s,t)}{\dd s} =  A^\top(s) P(s,t) + P(s,t) A(s) \\
~~~ + Q_2(s,t)  + C^\top(s) P(s,t) C(s)\\
~~~ - \bigl [ P(s,t) B_2(s)+ C^\top(s) P(s,t) D_2(s) \bigr ] \\
~~~ \times \bigl [ \widehat{R}_2(s,s) + D_2^\top(s) P(s,s)D_2(s) \bigr]^{-1} \\
~~~ \times \bigl [ B_2^\top(s) \widehat{P}(s,s) + D_2^\top(s) P(s,s) C(s) \bigr ],~ s \in [t,T] \\
- \frac{\dd Z(s,t)}{\dd s} = A^\top(s) Z(s,t) + Z(s,t) A(s) + \bar{Q}_2(s) \\
~~~ - Z(s,t) B_2(s) \bigl [ \widehat{R}_2(s,s) + D_2^\top(s) P(s,s)D_2(s) \bigr]^{-1} \\
~~~ \times \bigl [ B_2^\top(s) \widehat{P}(s,s) + D_2^\top(s) P(s,s) C(s) \bigr ],~ s \in [t,T]	 \\
P(T,t) = M_2(t),~ Z(T,t) = \bar{M}_2(t) \\
\widehat{P}(s,t) = P(s,t) + Z(s,t),
\end{cases}
\end{align}
which implies (see the notation defined in (\ref{eq_appendix_a_1}))
\begin{align}
\label{eq_18}
\begin{cases}
- \frac{\dd \widehat{P}(s,t)}{\dd s}	= A^\top(s) \widehat{P}(s,t) + \widehat{P}(s,t) A(s) \\
~~~ + \widehat{Q}_2(s,t)  + C^\top(s) P(s,t) + C(s) \\
~~~ - \bigl [ \widehat{P}(s,t) B_2(s) + C^\top(s) P(s,t) D_2(s) \bigr ] \\
~~~~~ \times \bigl [ \widehat{R}_2(s,s) + D_2^\top(s) P(s,s)D_2(s) \bigr]^{-1}  \\
~~~~~ \times \bigl [ B_2^\top(s) \widehat{P}(s,s) + D_2^\top(s) P(s,s) C(s) \bigr ] \\
\widehat{P}(T,t)  = \widehat{M}_2(t),~ s \in [t,T].
\end{cases}
\end{align}

Hence, by substituting (\ref{eq_15}) into the SDE in (\ref{eq_1}), with the notation defined in (\ref{eq_appendix_a_1}) of Appendix \ref{Appendix_A}, we have the following FBSDE that is equivalent to the FBSDE in (\ref{eq_11}) through the optimality condition in (\ref{eq_9}) and the transformation in (\ref{eq_12}) (note that $^*$ notation is dropped in $x$):
\begin{align}
\label{eq_20}
\begin{cases}
\dd x(s) = \Bigl [ H(s) x(s) + F(s)u(s) - G_1(s)  h(s) \\
~~~ -G_2 (s) l(s) \Bigr ] \dd s  + \Bigl [ \bar{H}(s) x(s) + \bar{F}(s) u(s) \\
~~~ - 	\bar{G}_1 (s) h(s) -  \bar{G}_2(s) l(s) \Bigr ] \dd B(s),~ s \in [t,T] \\
\dd h(s)  = \Bigl [ - \tilde{H}^\top(s,t) h(s) - \tilde{F}^\top(s,t) l(s)  \\
~~~ - K_1^\top(s,t) u(s) - K_2^\top(s,t) \mathbb{E}_t[u(s)]   \\
~~~ + \bar{K}_1^\top(s,t) \mathbb{E}_t[h(s)] + \bar{K}_2^\top(s,t) \mathbb{E}_t[l(s)] \Bigr ] \dd s  \\
~~~ + l(s) \dd B(s),~ s \in [t,T] \\
x(t)  = x_0,~ h(T)  = 0,
\end{cases}	
\end{align}
where the notation is defined in (\ref{eq_appendix_a_1}) of Appendix \ref{Appendix_A}. Note that $(h,l) \in \mathcal{L}_{\mathbb{F}}^2([0,T] \times [0,T],\mathbb{R}^n) \times \mathcal{L}_{\mathbb{F}}^2([0,T] \times [0,T],\mathbb{R}^n)$ is the BSDE. In fact, the FBSDE in (\ref{eq_20}) is the follower's rational behavior with respect to his equilibrium control under an arbitrary control of the leader. Note that $x$ in (\ref{eq_20}) depends on an arbitrary control of the leader, i.e., $x = x^* = x^{u,\alpha[u,x_0]}$, as in Definition \ref{Definition_2}(i). This is the optimization constraint of the leader's equilibrium control problem in Section \ref{Section_4}. 

In summary, we have the following result:

\begin{Theorem}\label{Theorem_1}
Consider \texttt{Problem LQ-FEC} with (\ref{eq_4}). Suppose that the coupled RDEs in (\ref{eq_17}) admit unique solutions, which are continuous in both variables. Then \texttt{Problem LQ-FEC} is solvable, and $v^*$ given in (\ref{eq_15}) is the state feedback representation of the (open-loop) equilibrium control of the follower, and $x$ in (\ref{eq_20}) is the corresponding equilibrium state process.
\end{Theorem}

\section{Leader's Equilibrium Control}\label{Section_4}
This section considers the leader's equilibrium control problem in the sense of Definition \ref{Definition_2}(ii). In the leader's problem, the objective functional is given in (\ref{eq_2}) with $v$ replaced by $v^*$ in (\ref{eq_15}) and the FBSDE in (\ref{eq_20}) is the corresponding optimization constraint (see Definition \ref{Definition_2}(ii)). As mentioned in Section \ref{Section_3}, (\ref{eq_20}) is the rational behavior of the follower with (\ref{eq_15}). The leader's problem is denoted by \texttt{Problem LQ-LEC} (LQ leader's equilibrium control problem). 

The following result states the leader's equilibrium control and the associated equilibrium state process in view of Definition \ref{Definition_2}(ii). 


\begin{Theorem}\label{Theorem_2}
Consider \texttt{Problem LQ-LEC} with (\ref{eq_4}). Assume that the assumptions in Theorem \ref{Theorem_1} hold. Suppose that the pair $(x^*,u^*) \in \mathcal{S}_\mathbb{F}([t,T],\mathbb{R}^n) \times \mathcal{U}[t,T]$ is the state-control pair of the leader with $x^*(t) = x_0$. Consider the following FBSDE:
\begin{align}
\label{eq_22}
\begin{cases}
\dd \phi(s,t) = \Bigl [ \tilde{H}(s,t) \phi(s,t) - \bar{K}_1(s,t) \mathbb{E}_t[\phi(s,t)] \\
~~~ + G_1^\top(s) y(s,t)  + \bar{G}_1^\top(s) z(s,t) \Bigr ] \dd s  \\
~~~ + \Bigl [ \tilde{F}(s,t) \phi(s) - \bar{K}_2(s,t) \mathbb{E}_t[\phi(s,t)]  + G_2^\top(s) y(s,t) \\
~~~ + \bar{G}_2^\top(s) z(s,t) \Bigr ] \dd B(s),~ s \in [t,T] \\
\dd y(s,t) = - \Bigl [ H^\top(s) y(s,t) + \bar{H}^\top(s) z(s,t) \\
~~~ + Q_1(s,t) x^*(s) + \bar{Q}_1(s,t) \mathbb{E}_t[x^*(s)] \Bigr ] \dd s \\
~~~ + z(s,t) \dd B(s),~ s \in [t,T] \\
\phi(t,t)  = 0 \\
y(T,t)  = M_1(t) x^*(T) + \bar{M}_1(t) \mathbb{E}_t[x^*(T)].
\end{cases}
\end{align}
Suppose that $(\phi,y,z) \in \mathcal{L}_{\mathbb{F}}^2([0,T] \times [0,T],\mathbb{R}^n) \times \mathcal{L}_{\mathbb{F}}^2([0,T] \times [0,T],\mathbb{R}^n) \times \mathcal{L}_{\mathbb{F}}^2([0,T] \times [0,T],\mathbb{R}^n)$ is the solution to the FBSDE in (\ref{eq_22}), which is continuous in both $s \in [t,T]$ and $t \in [0,T]$. Then the corresponding (open-loop) equilibrium control of the leader satisfies for $t \in [0,T]$,
\begin{align}
\label{eq_23}
& F^\top(t) y(t,t) + \bar{F}^\top(t) z(t,t) + \widehat{R}_1(t,t) u^*(t) = 0,
\end{align}
where $\widehat{R}_1(s,t) = R_1(s,t) + \bar{R}_1(s,t)$. In particular, for $t \in [0,T]$, we have
\begin{align}
\label{eq_24}
	u^*(t) & = - \widehat{R}_1^{-1}(t,t) ( F^\top(t) y(t,t) + \bar{F}^\top(t) z(t,t)).
\end{align}
\begin{Proof}
The proof is given in Appendix \ref{Appendix_B}.	
\end{Proof}
\end{Theorem}
Unlike the follower's case, due to the FBSDE constraint in the leader's problem, it is not clear that (\ref{eq_23}) is the necessary and sufficient condition for the leader's equilibrium control. Hence, the uniqueness of the leader's equilibrium control is not known yet.

Let (note that $^*$ notation is dropped in $x$) 
\begin{align*}
X(s) &= \begin{bmatrix}
 	x(s) \\
 	\phi(s,t)
 \end{bmatrix},~ Y(s,t) = \begin{bmatrix}
 	y(s,t) \\
 	h(s)
 \end{bmatrix},~L(s,t) = \begin{bmatrix}
 	z(s,t) \\
 	l(s)
 \end{bmatrix}.
\end{align*}
Note that $X(t) = \begin{bmatrix} x_0^\top & 0 \end{bmatrix}^\top$ in view of the initial condition in (\ref{eq_22}). Then with the notation defined in (\ref{eq_appendix_a_2}) of Appendix \ref{Appendix_A}, the FBSDEs in (\ref{eq_20}) and (\ref{eq_22}) can be written as 
\begin{align}
\label{eq_25}
\begin{cases}
\dd X(s) = \Bigl [ \mathcal{A}_1(s,t) X(s) + \mathcal{A}_2(s,t) \mathbb{E}_t[X(s)]  \\
~~~ + \mathcal{B}(s) u^*(s) + \mathcal{C}(s)Y(s,t) + \mathcal{D}(s) L(s,t) \Bigr ] \dd s \\
~~~ + \Bigl [ \bar{\mathcal{A}}_1(s,t) X(s) + \bar{\mathcal{A}}_2(s,t) \mathbb{E}_t[X(s)] + \bar{\mathcal{B}}(s) u(s) \\
~~~ + \bar{\mathcal{C}}(s)Y(s,t) + \bar{\mathcal{D}}(s) L(s,t) \Bigr ] \dd B(s),~ s \in [t,T] \\
\dd Y(s,t) = - \Bigl [ \mathcal{A}_1^\top(s,t)Y(s,t) + \bar{\mathcal{A}}_1^\top(s,t) L(s,t) \\
~~~ + \mathcal{A}_2^\top(s,t) \mathbb{E}_t[Y(s,t)]  + \bar{\mathcal{A}}_2^\top(s,t) \mathbb{E}_t[L(s,t)] \\
~~~ + \mathcal{Q}(s,t) X(s) + \bar{\mathcal{Q}}(s,t) \mathbb{E}_t[X(s)] \\
~~~ + \mathcal{G}^\top(s,t) u^*(s) + \bar{\mathcal{G}}^\top(s,t) \mathbb{E}_t[u^*(s)] \Bigr ] \dd s \\
~~~ + L(s,t) \dd B(s),~ s \in [t,T]	\\
X(t) = X_0 \\
 Y(T,t) = \begin{bmatrix}
	 M_1(t) x(T) + \bar{M}_1(t) \mathbb{E}_t[x(T)] \\
 	 0 
 \end{bmatrix},
\end{cases}
\end{align}
and the optimality condition in (\ref{eq_23}) is equivalent to 
\begin{align}
\label{eq_26}
 & \mathcal{B}^\top(t) Y(t,t) + \mathcal{\bar{B}}^\top(t) L(t,t) + \widehat{\mathcal{R}}(t,t) u^*(t) = 0,	
\end{align}
for $t \in [0,T]$, where $\widehat{\mathcal{R}}(s,t) = \widehat{R}_1(s,t) = R_1(s,t) + \bar{R}_1(s,t)$ (see the notation defined in (\ref{eq_appendix_a_2}) of Appendix \ref{Appendix_A}). 

We now obtain the state feedback representation of the equilibrium control of the leader via the generalized Four-Step Scheme and the coupled RDEs. Consider the following transformation:
\begin{align}
\label{eq_27}
Y(s,t) = \mathcal{P}(s,t)X(s)	+ \mathcal{Z}(s,t)\mathbb{E}_t[X(s)],~ s \in [t,T],
\end{align}
where $\mathcal{P}$ and $\mathcal{Z}$ are coupled RDEs with $\mathcal{P}, \mathcal{Z}:[0,T] \times [0,T] \rightarrow \mathbb{R}^{2n \times 2n}$, and their explicit expressions will be determined later. Unlike the follower's case in (\ref{eq_12}), there is no additional BSDE term in (\ref{eq_27}). Also, $\mathcal{P}$ and $\mathcal{Z}$ are $2n \times 2n$ dimensional due to the augmented state space $X$, $Y$ and $L$. From (\ref{eq_27}), we have $\mathcal{P}(T,t) = \mathcal{M}(t)$ and $\mathcal{Z}(T,t) = \bar{\mathcal{M}}(t)$. Let $\widehat{\mathcal{P}}(s,t) = \mathcal{P}(s,t) + \mathcal{Z}(s,t)$.

By using It\^o's formula, (note $\widehat{\mathcal{A}}_1(s,t) = \mathcal{A}_1(s,t) + \mathcal{A}_2(s,t)$)
\begin{align}
\label{eq_28}
& \dd Y(s,t)  = - \Bigl [ \mathcal{A}_1^\top(s,t)Y(s,t) + \bar{\mathcal{A}}_1^\top(s,t) L(s,t)    \\
&  + \mathcal{A}_2^\top(s,t) \mathbb{E}_t[Y(s,t)]  + \bar{\mathcal{A}}_2^\top(s,t) \mathbb{E}_t[L(s,t)] \nonumber \\
&  + \mathcal{Q}(s,t) X(s) + \bar{\mathcal{Q}}(s,t) \mathbb{E}_t[X(s)] \nonumber \\
&  + \mathcal{G}^\top(s,t) u^*(s) + \bar{\mathcal{G}}^\top(s,t) \mathbb{E}_t[u^*(s)] \Bigr ] \dd s + L(s,t) \dd B(s)	\nonumber \\
& = \frac{\dd \mathcal{P}(s,t)}{\dd s} X(s) \dd s + \frac{\dd \mathcal{Z}(s,t)}{\dd s} \mathbb{E}_t[X(s)] \dd s \nonumber \\
& + \mathcal{P}(s,t)\Bigl [ \mathcal{A}_1(s,t) X(s) + \mathcal{A}_2(s,t) \mathbb{E}_t[X(s)]  \nonumber \\
&  + \mathcal{B}(s) u^*(s) + \mathcal{C}(s)Y(s,t) + \mathcal{D}(s) L(s,t) \Bigr ] \dd s +  \mathcal{P}(s,t) \nonumber \\
&   \times \Bigl [ \bar{\mathcal{A}}_1(s,t) X(s) + \bar{\mathcal{A}}_2(s,t) \mathbb{E}_t[X(s)] + \bar{\mathcal{B}}(s) u^*(s) \nonumber \\
&  + \bar{\mathcal{C}}(s)Y(s,t) + \bar{\mathcal{D}}(s) L(s,t) \Bigr ] \dd B(s) \nonumber \\
& + \mathcal{Z}(s,t) \Bigl [ \widehat{\mathcal{A}}_1(s,t) \mathbb{E}_t[X(s)] + \mathcal{B}(s)\mathbb{E}_t[u^*(s)] \nonumber \\
& + \mathcal{C}(s) \mathbb{E}_t[Y(s,t)] + \mathcal{D}(s) \mathbb{E}_t[L(s,t)] \Bigr ] \dd s. \nonumber
\end{align}
Then we can easily see that $L(s,t) = \mathcal{P}(s,t) \bigl [ \bar{\mathcal{A}}_1(s,t) X(s) + \bar{\mathcal{A}}_2(s,t) \mathbb{E}_t[X(s)]  + \bar{\mathcal{B}}(s) u^*(s)  + \bar{\mathcal{C}}(s) \mathcal{P}(s,t) X(s) + \bar{\mathcal{C}}(s) \mathcal{Z}(s,t) \mathbb{E}_t[X(s)] + \bar{\mathcal{D}}(s) L(s,t) \bigr ]$,
which, together with the existence of $ [ I - \mathcal{P}(s,t) \bar{\mathcal{D}}(s) ]^{-1}$, implies
\begin{align}
\label{eq_29}
L(s,t) &= \bigl [ I - \mathcal{P}(s,t) \bar{\mathcal{D}}(s) \bigr ]^{-1} \mathcal{P}(s,t) \\
&~~~ \times 	 \bigl [ \bar{\mathcal{A}}_1(s,t) X(s) + \bar{\mathcal{A}}_2(s,t) \mathbb{E}_t[X(s)] + \bar{\mathcal{B}}(s) u^*(s) \nonumber \\
&~~~   + \bar{\mathcal{C}}(s) \mathcal{P}(s,t) X(s) + \bar{\mathcal{C}}(s) \mathcal{Z}(s,t) \mathbb{E}_t[X(s)] \bigr ].\nonumber
\end{align}

By substituting (\ref{eq_27}) and (\ref{eq_29}) into the optimality condition in (\ref{eq_26}) (note $\widehat{\mathcal{P}}(s,t) = \mathcal{P}(s,t) + \mathcal{Z}(s,t)$ and $\widehat{\mathcal{A}}_2(s,t) = \bar{\mathcal{A}}_1(s,t) + \bar{\mathcal{A}}_2(s,t)$; see also the notation defined in (\ref{eq_appendix_a_2}))
\begin{align*}
 & \mathcal{B}^\top(t) \widehat{\mathcal{P}}(t,t)X(t) +  \widehat{\mathcal{R}}(t,t) u^*(t) \\
 & + \bar{\mathcal{B}}^\top(t) \bigl [ I - \mathcal{P}(t,t) \bar{\mathcal{D}}(t) \bigr ]^{-1}  \mathcal{P}(t,t)\\
& \times 	 \bigl [ \widehat{\mathcal{A}}_2(t,t) X(t)   + \bar{\mathcal{B}}(t) u^*(t)  + \bar{\mathcal{C}}(t) \widehat{\mathcal{P}}(t,t) X(t) \bigr ] = 0,
\end{align*}
which implies
\begin{align}
\label{eq_30}
u^*(t) &= - \bigl [ \widehat{\mathcal{R}}(t,t) + \bar{\mathcal{B}}^\top(t) \bigl [ I - \mathcal{P}(t,t) \bar{\mathcal{D}}(t) \bigr ]^{-1}  \mathcal{P}(t,t)  \bar{\mathcal{B}}(t) \bigr ]^{-1} \nonumber
\\
&~~~ \times \bigl [ \mathcal{B}^\top(t) \widehat{\mathcal{P}}(t,t) + \bar{\mathcal{B}}^\top(t) \bigl [ I - \mathcal{P}(t,t) \bar{\mathcal{D}}(t) \bigr ]^{-1}  \nonumber \\
&~~~~~ \times \mathcal{P}(t,t) (\widehat{\mathcal{A}}_2(t,t) + \bar{\mathcal{C}}(t) \widehat{\mathcal{P}}(t,t) )\bigr ] X(t),
\end{align}
where for $t \in [0,T]$, we denote 
\begin{align}
\label{eq_31}
u^*(t) &= - \Pi(t,t) X(t) = -\bar{\Pi}(t,t) x(t).	
\end{align}
Note that $\Pi$ is $m_1 \times 2n$, and $\bar{\Pi}$ is $m_1 \times n$. The last equality of (\ref{eq_31}) follows from the definition of $X$ and its initial condition. 
 Note that (\ref{eq_30}) (equivalently (\ref{eq_31})) is the state feedback representation of the leader's equilibrium control in (\ref{eq_24}) that is obtained by the transformations in (\ref{eq_27}) and (\ref{eq_29}).

From (\ref{eq_31}) and (\ref{eq_29}), we have
\begin{align}
\label{eq_32}
L(s,t) &= \bigl [ I - \mathcal{P}(s,t) \bar{\mathcal{D}}(s) \bigr ]^{-1} \mathcal{P}(s,t) \bigl [ \bar{\mathcal{A}}_1(s,t) X(s)    \\
&~~~ + \bar{\mathcal{A}}_2(s,t) \mathbb{E}_t[X(s)] - \bar{\mathcal{B}}(s) \Pi(s,s) X(s) \nonumber \\
&~~~   + \bar{\mathcal{C}}(s) \mathcal{P}(s,t) X(s) \nonumber + \bar{\mathcal{C}}(s) \mathcal{Z}(s,t) \mathbb{E}_t[X(s)] \bigr ]. \nonumber 
\end{align}
Substituting (\ref{eq_27}), (\ref{eq_31}) and (\ref{eq_32}) into $X$ in (\ref{eq_25}) yields
\begin{align}
\label{eq_33_1}
\begin{cases}
\dd X(s) = \Bigl [ \mathcal{A}_1(s,t) X(s)  + \mathcal{A}_2(s,t) \mathbb{E}_t[X(s)]  \\
~~~ - \mathcal{B}(s) \Pi(s,s) X(s)  + \mathcal{C}(s)\mathcal{P}(s,t)X(s)\\
 ~~~ 	+ \mathcal{C}(s) \mathcal{Z}(s,t)\mathbb{E}_t[X(s)]  \\
~~~ + \mathcal{D}(s) \bigl [ I - \mathcal{P}(s,t) \bar{\mathcal{D}}(s) \bigr ]^{-1}  \\
~~~ \times 	\mathcal{P}(s,t) \bigl [ \bar{\mathcal{A}}_1(s,t) X(s) + \bar{\mathcal{A}}_2(s,t) \mathbb{E}_t[X(s)]  \\
~~~ - \bar{\mathcal{B}}(s) \Pi(s,s) X(s)  + \bar{\mathcal{C}}(s) \mathcal{P}(s,t) X(s)  \\
~~~ + \bar{\mathcal{C}}(s) \mathcal{Z}(s,t) \mathbb{E}_t[X(s)] \bigr ]  \Bigr ] \dd s	 \\
~~~ +  \Bigl [ \bar{\mathcal{A}}_1(s,t) X(s) + \bar{\mathcal{A}}_2(s,t) \mathbb{E}_t[X(s)]  \\
 ~~~ - \bar{\mathcal{B}}(s) \Pi(s,s) X(s) + \bar{\mathcal{C}}(s) \mathcal{P}(s,t) X(s)   \\
~~~ + \bar{\mathcal{C}}(s) \mathcal{Z}(s,t) \mathbb{E}_t[X(s)] \\
~~~ + \bar{\mathcal{D}}(s) \bigl [ I - \mathcal{P}(s,t) \bar{\mathcal{D}}(s) \bigr ]^{-1} \\
~~~ \times 	\mathcal{P}(s,t) \bigl [ \bar{\mathcal{A}}_1(s,t) X(s) + \bar{\mathcal{A}}_2(s,t) \mathbb{E}_t[X(s)]  \\
~~~ - \bar{\mathcal{B}}(s) \Pi(s,s) X(s)  + \bar{\mathcal{C}}(s) \mathcal{P}(s,t) X(s)  \\
~~~ + \bar{\mathcal{C}}(s) \mathcal{Z}(s,t) \mathbb{E}_t[X(s)] \bigr ] \Bigr ] \dd B(s),~ s \in [t,T] \\
X(t) = X_0.
\end{cases}
\end{align}
By substituting (\ref{eq_27}), (\ref{eq_31}) and (\ref{eq_32}) into (\ref{eq_28}), we have
\begin{align*}
& \dd Y(s,t) = - \Bigl [ \mathcal{A}_1^\top(s,t)\mathcal{P}(s,t)X(s)	+ \mathcal{A}_1^\top(s,t) \\
& \times \mathcal{Z}(s,t)\mathbb{E}_t[X(s)]  + \bar{\mathcal{A}}_1^\top(s,t) \bigl [ I - \mathcal{P}(s,t) \bar{\mathcal{D}}(s) \bigr ]^{-1} \\
& \times 	\mathcal{P}(s,t) \Bigl [ \bar{\mathcal{A}}_1(s,t) X(s) + \bar{\mathcal{A}}_2(s,t) \mathbb{E}_t[X(s)] \\
& - \bar{\mathcal{B}}(s) \Pi(s,s) X(s)  + \bar{\mathcal{C}}(s) \mathcal{P}(s,t) X(s) \\
& + \bar{\mathcal{C}}(s) \mathcal{Z}(s,t) \mathbb{E}_t[X(s)] \Bigr ]  + \mathcal{A}_2^\top(s,t) \widehat{\mathcal{P}}(s,t) \mathbb{E}_t[X(s,t)]  \\
& + \bar{\mathcal{A}}_2^\top(s,t) \bigl [ I - \mathcal{P}(s,t) \bar{\mathcal{D}}(s) \bigr ]^{-1} \mathcal{P}(s,t) \Bigl [ \widehat{\mathcal{A}}_2(s,t) \\
& + \bar{\mathcal{C}}(s) \widehat{\mathcal{P}}(s,t)  	 - \bar{\mathcal{B}}(s) \Pi(s,s)    \Bigr ]	 \mathbb{E}_t[X(s)] + \mathcal{Q}(s,t) X(s)  \\
&  + \bar{\mathcal{Q}}(s,t) \mathbb{E}_t[X(s)] - \mathcal{G}^\top(s,t) \Pi(s,s)    \\
&  - \bar{\mathcal{G}}^\top(s,t) \Pi(s,s) \mathbb{E}_t[X(s)] \Bigr ] \dd s  + L(s,t) \dd B(s),
\end{align*}
and
\begin{align*}
\dd Y(s,t) & = \frac{\dd \mathcal{P}(s,t)}{\dd s} X(s) \dd s + \frac{\dd \mathcal{Z}(s,t)}{\dd s} \mathbb{E}_t[X(s)] \dd s \\
&~~~ + \mathcal{P}(s,t) \dd X(s) +\mathcal{Z}(s,t) \dd \mathbb{E}_t[X(s)], \nonumber 
\end{align*}
where the expressions of $\dd X(s)$ and $\dd \mathbb{E}_t[X(s)]$ can be obtained from (\ref{eq_33_1}). 

We compare the coefficients in the above two equalities. Then we can easily show that $\mathcal{P}$ and $\mathcal{Z}$ have to satisfy the following nonsymmetric coupled RDEs:
\begin{align}
\label{eq_34}
\begin{cases}
- \frac{\dd \mathcal{P}(s,t)}{\dd s} = \Lambda_1(s,t, \mathcal{P}(s,t), \mathcal{Z}(s,t)),~ s \in [t,T]\\
- \frac{\dd \mathcal{Z}(s,t)}{\dd s} = \Lambda_2(s,t, \mathcal{P}(s,t), \mathcal{Z}(s,t)),~ s \in [t,T]\\
\mathcal{P}(T,t) = \mathcal{M}(t),~ \mathcal{Z}(T,t) = \bar{\mathcal{M}}(t)		\\
	\det (I - \mathcal{P}(s,t) \bar{\mathcal{D}}(s))  \neq 0 \\
	\det \bigl (\widehat{\mathcal{R}}(t,t)  + \bar{\mathcal{B}}^\top(t) \bigl [ I - \mathcal{P}(t,t) \bar{\mathcal{D}}(t) \bigr ]^{-1}  \mathcal{P}(t,t)  \bar{\mathcal{B}}(t) \bigr) \neq 0\\
		\widehat{\mathcal{P}}(s,t) = \mathcal{P}(s,t) + \mathcal{Z}(s,t),
\end{cases}	
\end{align}
where the explicit expressions of $\Lambda_1$ and $\Lambda_2$ are provided in Appendix \ref{Appendix_C}. In (\ref{eq_34}), $\widehat{\mathcal{P}}(s,t)= \mathcal{P}(s,t) + \mathcal{Z}(s,t)$ satisfies $-\frac{\dd \widehat{\mathcal{P}}(s,t)}{\dd s} = \Lambda_3(s,t,\mathcal{P}(s,t),\mathcal{Z}(s,t))$, $s \in [t,T]$, with $\widehat{\mathcal{P}}(T,t) = \mathcal{M}(t) + \bar{\mathcal{M}}(t)$ given (\ref{eq_appendix_c_3}) in Appendix \ref{Appendix_C}. Note that in (\ref{eq_34}), the two nonsingularity conditions are included due to the invertibility of the matrices in (\ref{eq_32}) and (\ref{eq_30}).


By substituting $u^*$ in (\ref{eq_31}) (equivalently (\ref{eq_30})) into the leader's optimization constraint $x$ in (\ref{eq_20})\footnote{As mentioned at the beginning of Section \ref{Section_4}, (\ref{eq_20}) is the rational behavior of the follower with respect to his equilibrium control in (\ref{eq_15}).} (equivalently (\ref{eq_25})), we have (note that $^*$ notation is dropped in $x$):
\begin{align}
\label{eq_35}
\begin{cases}
\dd x(s) = \Bigl [ \Psi(s) x(s)  - G_1(s) h(s)  - G_2(s) l(s) \Bigr ]\dd s  \\
~~~ + \Bigl [ \bar{\Psi}(s) x(s) - \bar{G}_1(s) h(s) - \bar{G}_2(s) l(s) \Bigr ] \dd B(s) \\
x(t)  = x_0,~s \in [t,T],
\end{cases}
\end{align}
where $\Psi(s) = H(s) - F(s) \bar{\Pi}(s,s)$, $\bar{\Psi}(s) = \bar{H}(s) - \bar{F}(s) \bar{\Pi}(s,s)$, and $(h,l)$ is given in (\ref{eq_20}) (see also (\ref{eq_25})). Note that (\ref{eq_35}) is the leader's equilibrium state process, where (\ref{eq_31}) (equivalently in (\ref{eq_30})) is the corresponding the leader's equilibrium control in view of Definition \ref{Definition_2}(ii).

In summary, we have the following result:
\begin{Theorem}\label{Theorem_3}
Consider \texttt{Problem LQ-LEC} with (\ref{eq_4}). Suppose that the assumptions in Theorem \ref{Theorem_1} hold. Assume that the coupled RDEs in (\ref{eq_34}) admit unique solutions, which are continuous in both variables. Then \texttt{Problem LQ-LEC} is solvable, and 
$u^* \in \mathcal{U}[t,T]$ given in (\ref{eq_31}) (equivalently in (\ref{eq_30})) is the state feedback representation of the (open-loop) equilibrium control of the leader, where $x$ in (\ref{eq_35}) is the corresponding equilibrium state process.
\end{Theorem}

In view of Theorems \ref{Theorem_1} and \ref{Theorem_3}, we now state the existence of the time-consistent Stackelberg equilibrium of \texttt{Problem LQ-TI-MF-SDG} (see Definition \ref{Definition_2}(iii)).

\begin{Corollary}\label{Corollary_1}
	Consider \texttt{Problem LQ-TI-MF-SDG} of the paper with (\ref{eq_4}). Assume that the assumptions in Theorems \ref{Theorem_1} and \ref{Theorem_3} hold. Let $u^*$ be given in (\ref{eq_31}) (equivalently (\ref{eq_30})), and $v^*$ given in (\ref{eq_15}) with $u$ replaced by $u^*$ in (\ref{eq_31}), i.e., $v^* = \alpha[u^*,x_0] \in \mathcal{V}[t,T]$. Then $(u^*,v^*) \in \mathcal{U}[t,T] \times \mathcal{V}[t,T]$ constitutes the (adapted open-loop) time-consistent Stackelberg equilibrium of \texttt{Problem LQ-TI-MF-SDG}, and (\ref{eq_35}) is the corresponding equilibrium state process.
\end{Corollary}

	In view of Theorems \ref{Theorem_1} and \ref{Theorem_3} (and Corollary \ref{Corollary_1}), the solvability (existence and uniqueness of the solution) of the nonsymmetric coupled RDEs of the follower (\ref{eq_17}) and the leader (\ref{eq_34}) is crucial to characterize the time-consistent Stackelberg equilibrium. However, their general solvability problem requires a different approach than that for the RDEs in various classes of (time-consistent) LQ optimal control and differential games, since (\ref{eq_17}) and (\ref{eq_34}) have two time variables and the cost parameters could be general nonexponential discounting depending on the initial time. Such coupled RDEs have not been studied in the existing literature, which we will address in the future. Note that the solvability of (\ref{eq_17}) and (\ref{eq_34}) is verified numerically in Section \ref{Section_5} for the time-inconsistent resource-allocation Stackelberg game.

\section{Numerical Examples}\label{Section_5}

\begin{figure}[!t]
\centering
\includegraphics[scale=0.43]{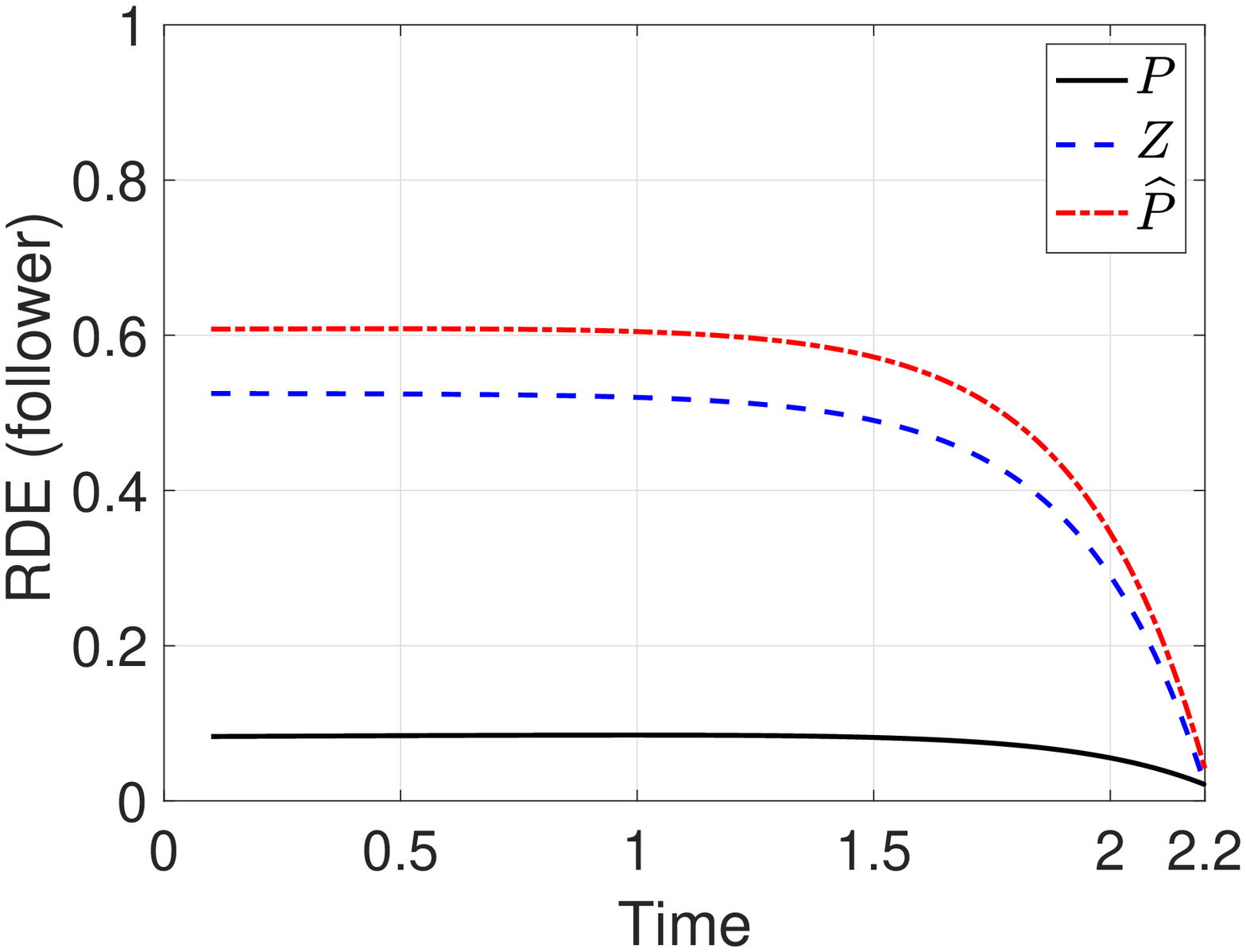}
\includegraphics[scale=0.43]{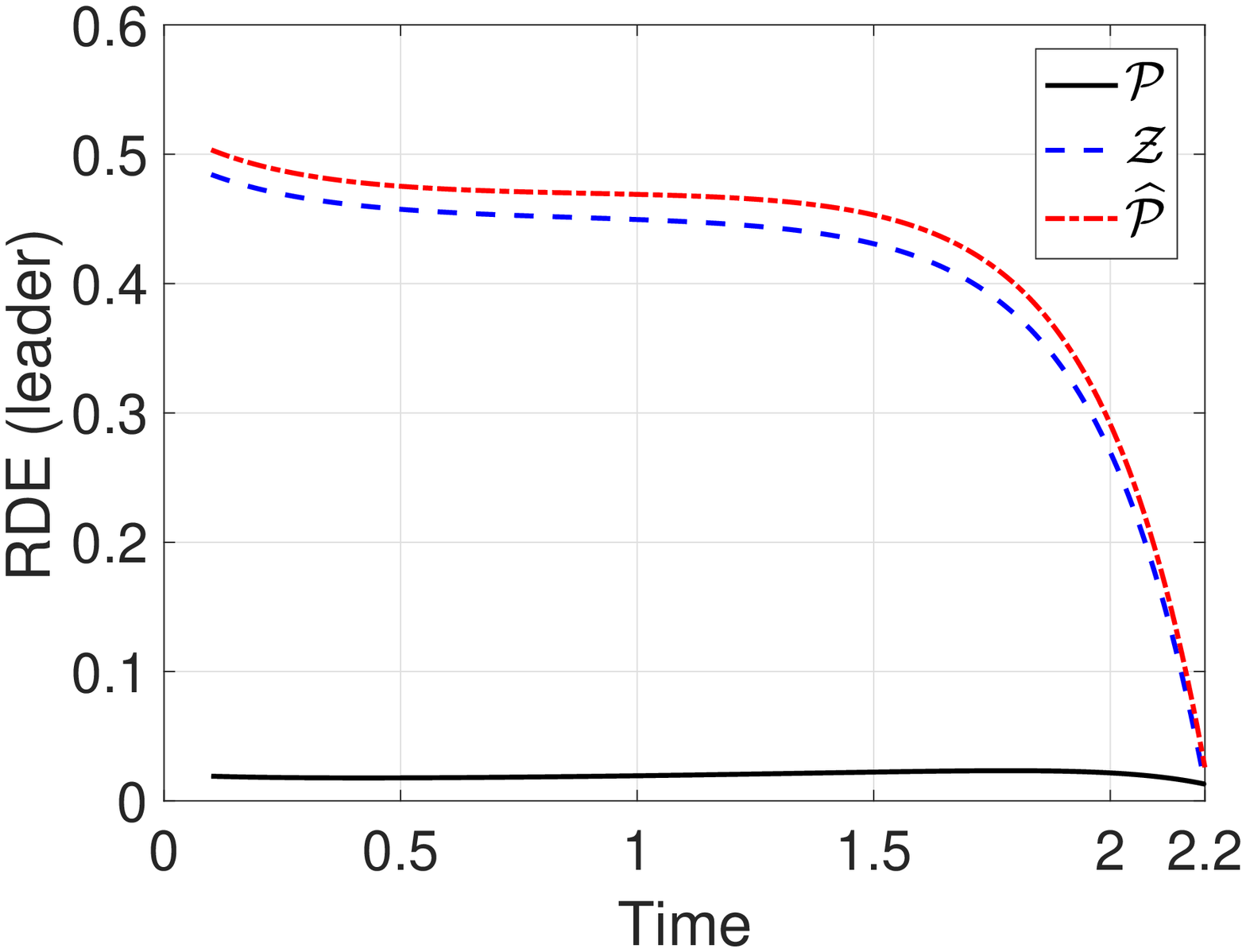}
\caption{The evolutions of the nonsymmetric coupled RDEs of the follower and the leader for Case I (top: follower, bottom: leader).}
\label{Fig_1}
\end{figure}

\begin{figure}[!t]
\centering
\includegraphics[scale=0.43]{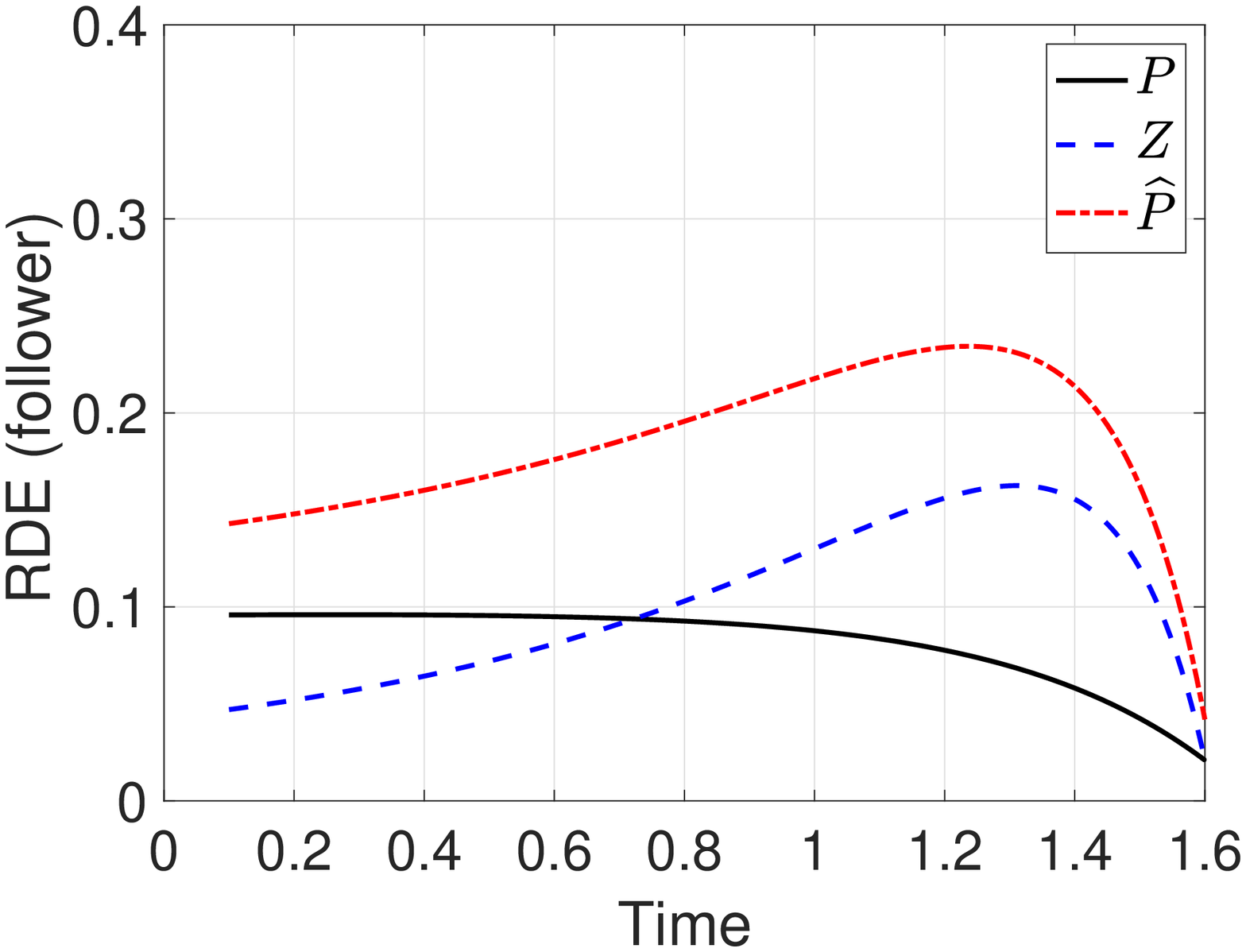}
\includegraphics[scale=0.43]{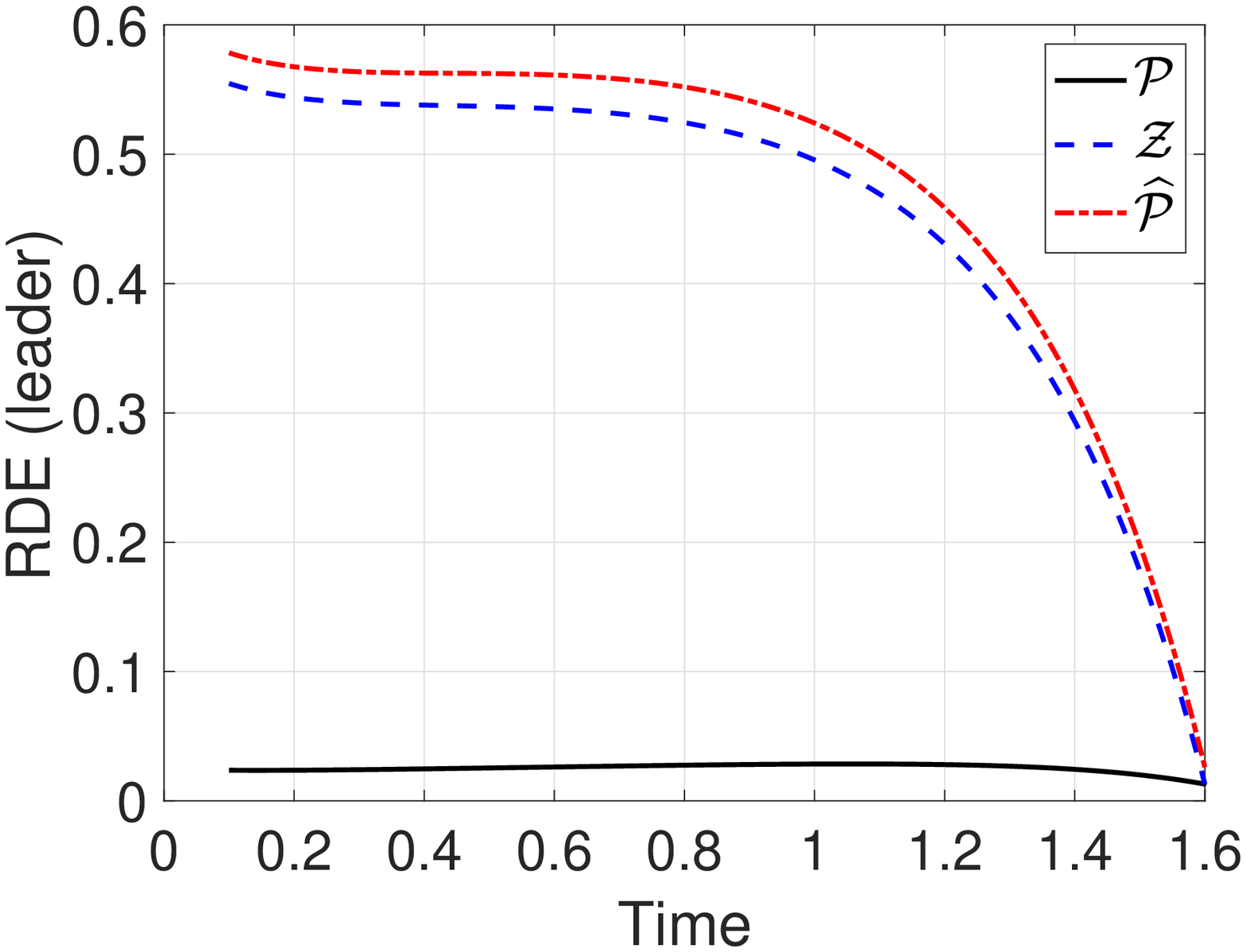}
\caption{The evolutions of the nonsymmetric coupled RDEs of the follower and the leader for Case II (top: follower, bottom: leader).}
\label{Fig_2}
\end{figure}

As mentioned, it is necessary to check the solvability of the nonsymmetric coupled RDEs of the follower (\ref{eq_17}) and the leader (\ref{eq_34}) to characterize the time-consistent Stackelberg equilibrium. This section provides the numerical examples to check their solvability.


We consider the modified resource-allocation Stackelberg game considered in \cite{Chen_TAC_1972,Lin_TAC_2018}, where unlike \cite{Lin_TAC_2018} the time-inconsistent problem is formulated. Specifically, the conditional expectations (mean-field variables) are included in the objective functionals, and the nonexponential discounting is hyperbolic similar to  \cite{Ekeland_SIFM_2012}. We consider two different cases:
\begin{itemize}
\setlength{\itemindent}{2.0em}
\item[Case I] $t=0.1$, $T=2.2$, $A=-1.6$, $B_1=-0.3$, $B_2 = 1$, $C=0.5$, $D_1 = D_2 = 0$, $\bar{Q}_1 = \bar{Q}_2 = 2 $, $\bar{R}_1 = \bar{R}_2 = 0$, $R_1(s,t) = (1+0.7(s-t))^{-1.8}$, $R_2(s,t) = (1+0.2(s-t))^{-0.3}$, $Q_1(s,t) = (5+1.2(s-t))^{-1.2}$,  $Q_2(s,t) = (10+0.8(s-t))^{-0.5}$, $M_1(t) = \bar{M}_1(t) = 1.3 t^2$, $M_2(t) = \bar{M}_2(t) = 2.1t^2$;
\item[Case II] $t=0.1$, $T=1.6$, $A=-1.6$, $B_1=-0.3$, $B_2 = 1$, $C=0.5$, $D_1 = 0.7$, $D_2 = 0.2$, $Q_1(s,t) = (5+1.2(s-t))^{-1.2}$,  $R_1(s,t) = (1+0.7(s-t))^{-1.8}$, $\bar{Q}_1 = 2 $, $\bar{R}_1(s,t) = (20+3.7(s-t))^{-0.1}$, $Q_2(s,t) = (10+0.8(s-t))^{-0.5}$, $R_2(s,t) = (1+0.2(s-t))^{-0.3}$, $\bar{Q}_2 = (0.7 + 3.1(s-t))^{-1.3}$, $\bar{R}_2 = 25$, $M_1(t) = \bar{M}_1(t) = 1.3 t^2$, $M_2(t) = \bar{M}_2(t) = 2.1t^2$.
\end{itemize}
Note the hyperbolic cost parameters, $R_i$, $\bar{R}_i$, $Q_i$, $\bar{Q}_i$, $M_i$ and $\bar{M}_i$, in Case I and Case II. It can be seen that Case II is more general in that we consider the controlled diffusion term due to nonzero $D_1$ and $D_2$, and that the inclusion of the conditional expectation of the controls in the objective functionals. Note that with Case I and Case II, (\ref{eq_4}) holds. Using the Euler's method (see \cite{Kahaner_numerical_book}), we obtain the plots of the evolutions of the nonsymmetric coupled RDEs of (\ref{eq_17}) and (\ref{eq_34}) for Cases I and II, which are depicted in Figs. \ref{Fig_1} and \ref{Fig_2}. These simulation results show that the solutions of the coupled RDEs are well defined and do not have a finite escape time \cite{Liberzon_Book}. This implies that in view of Theorems \ref{Theorem_1} and \ref{Theorem_3}, and Corollary \ref{Corollary_1}, the time-consistent Stackelberg equilibrium for the time-inconsistent resource-allocation problem with Case I and Case II exists, which can be characterized by substituting the parameters in the corresponding case into $u^*$ in (\ref{eq_31}) (equivalently (\ref{eq_30})) and $v^*$ in (\ref{eq_15}) with $u$ replaced by $u^*$ in (\ref{eq_31}).

\section{Conclusions}\label{Section_6}
In this paper, we have considered the linear-quadratic (LQ) time-inconsistent mean-field Stackelberg differential game (\texttt{Problem LQ-TI-MF-SDG}). The (time-consistent) equilibrium control of the follower and its state feedback representation in terms of the nonsymmetric coupled RDEs and the BSDE have been obtained. The equilibrium state process of the follower is the FBSDE that captures his rational behavior. Then the leader's (time-consistent) equilibrium control has been obtained via the variational method under the constraint of the follower's rational behavior (the FBSDE). The state feedback representation of the leader's equilibrium control in terms of the nonsymmetric coupled RDEs has been characterized via the generalized decoupling method. Under the solvability of the nonsymmetric coupled RDEs of the leader and the follower, the equilibrium controls of the leader and the follower constitute the time-consistent Stackelberg equilibrium of \texttt{Problem LQ-TI-MF-SDG}. A future research topic is to consider the discrete-time time-inconsistent leader-follower Stackelberg game that can be viewed as an extension of \cite{Ni_TAC_2018}. Another topic would be the existence and uniqueness of the solution to (\ref{eq_17}) and (\ref{eq_34}) and their numerical computation approaches. Note that under some assumptions on the coefficients in (\ref{eq_17}) and (\ref{eq_34}), the Picard fixed point argument can be used to obtain both existence and uniqueness.

\newpage
\appendices
\numberwithin{equation}{section}
\renewcommand{\theequation}{\thesection.\arabic{equation}}

\section{Notations for Sections \ref{Section_3} and \ref{Section_4}}\label{Appendix_A}

The following notation is used in Section \ref{Section_3}:
\begin{equation*}
\begin{aligned}
S(s,s)  &= [ \widehat{R}_2(s,s) + D_2^\top(s) P(s,s)D_2(t) \bigr] \\
H(s) &= A(s) - B_2(s) S^{-1}(s,s) \\
&~~~ \times \bigl [ B_2^\top(s) \widehat{P}(s,s) + D_2^\top(t) P(s,s) C(s) \bigr ]  \\
\bar{H}(s) &= C(s) - D_2(s) S^{-1}(s,s) \\
&~~~ \times \bigl [ B_2^\top(s) \widehat{P}(s,s) + D_2^\top(t) P(s,s) C(s) \bigr ]  \\
\tilde{H}(s,t) &= A(s) - B_2(s) (S^{-1}(s,s))^\top \\
&~~~ \times \bigl [ B_2^\top(s) P^\top(s,t)  + D_2^\top(s) P^\top(s,t) C(s) \bigr ]   \\
F(s) &= B_1(s) - B_2(s) S^{-1}(s,s) D_2^\top (s) P(s,s) D_1(s) \\
\bar{F}(s) &=D_1(s) - D_2(s) S^{-1}(s,s) D_2^\top (s) P(s,s) D_1(s) \\
\tilde{F}(s,t) &= C(s) - D_2(s) (S^{-1}(s,s))^\top \\
&~~~ \times \bigl [ B_2^\top(s) P^\top(s,t)  + D_2^\top(s) P^\top(s,t) C(s) \bigr ] \\
G_1(s) &= B_2(s) S^{-1}(s,s)  B_2^\top(s) \\
G_2(s) &= B_2(s) S^{-1}(s,s)  D_2^\top(s) = \bar{G}_1^\top(s) \\
\bar{G}_2(s) &= D_2(s) S^{-1}(s,s)  D_2^\top(s) \\
K_1(s,t)  &= D_1^\top(s) P^\top(s,t) C(s) + B_1^\top(s) P^\top(s,t) \\
& ~~~ - D_1^\top(s) P^\top(s,s) D_2(s) (S^{-1}(s,s))^\top \\
& ~~~   \times \bigl [ B_2^\top(s) P^\top(s,t)  + D_2^\top(s) P^\top(s,t) C(s) \bigr ]
\end{aligned}
\end{equation*}
\begin{equation}
\label{eq_appendix_a_1}	
\begin{aligned}
K_2(s,t)  &= B_1^\top(s) Z^\top(s,t) - D_1^\top(s) P^\top(s,s) D_2(s) \\
& ~~~ \times (S^{-1}(s,s))^\top B_2^\top(s) Z^\top(s,t) \\
\bar{K}_1(s,t)  &= B_2(s) (S^{-1}(s,s))^\top(s,s)  B_2^\top(s) Z^\top(s,t) \\
\bar{K}_2(s,t)  &= D_2(s) (S^{-1}(s,s))^\top   B_2^\top(s) Z^\top(s,t) \\
\widehat{R}_2(s,t) &= R_2(s,t) + \bar{R}_2(s,t) \\
\widehat{M}_2(t) &= M_2(t) + \bar{M}_2(t) \\
\widehat{Q}_2(s,t) &= Q_2(s,t) + \bar{Q}_2(s,t)\\
\widehat{P}(s,t) &= P(s,t) + Z(s,t).
\end{aligned}
\end{equation}

The following notation is used in Section \ref{Section_4}:
\begin{equation}
\label{eq_appendix_a_2}
\begin{aligned}
\mathcal{A}_1(s,t)  &= \begin{bmatrix}
 H(s) & 0 \\
 0 & \tilde{H}(s,t)	
 \end{bmatrix} \\
 \bar{\mathcal{A}}_1(s,t)  &= \begin{bmatrix}
 \bar{H}(s) & 0 \\
 0 & \tilde{F}(s,t)	
 \end{bmatrix}
\\
\mathcal{A}_2(s,t)  &= \begin{bmatrix}
 0 & 0 \\
 0 & -\bar{K}_1(s,t)
 \end{bmatrix} \\
 \bar{\mathcal{A}}_2(s,t) &= \begin{bmatrix}
 0 & 0 \\
 0 & -\bar{K}_2(s,t)	
 \end{bmatrix}
\\
\mathcal{B}(s) &= \begin{bmatrix}
 F(s) \\
 0	
 \end{bmatrix},~ \bar{\mathcal{B}}(s) = \begin{bmatrix}
 \bar{F}(s) \\
 0	
 \end{bmatrix} \\
\mathcal{C}(s)  &= \begin{bmatrix}
 0 & -G_1(s) \\
 G_1^\top(s) & 0	
 \end{bmatrix} \\
 \bar{\mathcal{C}}(s)  &= \begin{bmatrix}
 	0 & -\bar{G}_1(s) \\
 	G_2^\top(s) & 0
 \end{bmatrix} \\
\mathcal{D}(s) &= \begin{bmatrix}
 	0 & -G_2(s) \\
 	\bar{G}_1^\top(s) & 0
 \end{bmatrix} \\
 \bar{\mathcal{D}}(s) &= \begin{bmatrix}
 	0 & -\bar{G}_2(s) \\
 	\bar{G}_2^\top(s) & 0
 \end{bmatrix} \\
 \mathcal{Q}(s,t) & = \begin{bmatrix}
 	Q_1(s,t) & 0 \\
 	0 & 0
 \end{bmatrix},~ \bar{\mathcal{Q}}(s,t) = \begin{bmatrix}
 	\bar{Q}_1(s,t) & 0 \\
 	0 & 0
 \end{bmatrix} \\
 \mathcal{G}(s,t)  &= \begin{bmatrix}
 	0 &
 	K_1(s,t)
 \end{bmatrix},~  \bar{\mathcal{G}}(s,t) = \begin{bmatrix}
 	0 &
 	K_2(s,t)
 \end{bmatrix} \\
 \mathcal{M}(s) &= \begin{bmatrix}
 	M_1(s) & 0 \\
 	0 & 0
 \end{bmatrix},~ \bar{\mathcal{M}}(s) = \begin{bmatrix}
 	\bar{M}_1(s) & 0 \\
 	0 & 0
 \end{bmatrix} \\
 X_0  &= \begin{bmatrix}
 	x_0 \\
 	0
 \end{bmatrix} \\
 \widehat{\mathcal{R}}(s,t) &= \widehat{R}_1(s,t) = R_1(s,t) + \bar{R}_1(s,t) \\
  \widehat{\mathcal{A}}_1(s,t) &= \mathcal{A}_1(s,t) + \mathcal{A}_2(s,t) \\
 \widehat{\mathcal{A}}_2(s,t) &= \bar{\mathcal{A}}_1(s,t) + \bar{\mathcal{A}}_2(s,t).
 \end{aligned}
\end{equation}

\section{Proof of Theorem \ref{Theorem_2}}\label{Appendix_B}
This appendix provides the proof of Theorem \ref{Theorem_2}.

\textit{Proof of Theorem \ref{Theorem_2}:} Let $u^*$ be given in (\ref{eq_24}), and the triplet $(x^*,h^*,l^*)$ be the FBSDE in (\ref{eq_20}) generated by $u^*$ (that can be obtained by substituting  (\ref{eq_24}) into (\ref{eq_20})).\footnote{Note that in (\ref{eq_20}), $^*$ notation is dropped in $x$.} For any $s \in [t,T)$, let $u^\epsilon$ be the control defined in (\ref{eq_6}), and the triplet $(x^\epsilon,h^\epsilon,l^\epsilon)$ be the FBSDE in (\ref{eq_20}) generated by $u^\epsilon$, i.e.,
\begin{align*}
\begin{cases}
\dd x^\epsilon(s) = \Bigl [ H(s) x^\epsilon(s) + F(s)u^\epsilon(s) - G_1(s)  h^\epsilon(s) \\
~~~ -G_2 (s) l^\epsilon(s) \Bigr ] \dd s  + \Bigl [ \bar{H}(s) x^\epsilon(s) + \bar{F}(s) u^\epsilon(s) \\
~~~ - 	\bar{G}_1 (s) h^\epsilon(s) -  \bar{G}_2(s) l^\epsilon(s) \Bigr ] \dd B(s),~ s \in [t,T] \\
\dd h^\epsilon(s)  = \Bigl [ - \tilde{H}^\top(s,t) h^\epsilon(s) - \tilde{F}^\top(s,t) l^\epsilon(s)  \\
~~~ - K_1^\top(s,t) u^\epsilon(s) - K_2^\top(s,t) \mathbb{E}_t[u^\epsilon(s)]   \\
~~~ + \bar{K}_1^\top(s,t) \mathbb{E}_t[h^\epsilon(s)] + \bar{K}_2^\top(s,t) \mathbb{E}_t[l^\epsilon(s)] \Bigr ] \dd s  \\
~~~ + l^\epsilon(s) \dd B(s),~ s \in [t,T] \\
x^\epsilon(t)  = x_0,~ h(T)  = 0.
\end{cases}	
\end{align*}
Let $\delta x = x^* - x^\epsilon$ $\delta h = h^* - h^\epsilon$, and $\delta l  = l^* - l^\epsilon$. Note that $\delta x(t) = 0$ and $\delta h(T) = 0$. 

By using It\^o's formula, we have
\begin{align*}
& \dd \langle \delta x(s),y(s,t) \rangle \\
& = \langle F(s) (u^\epsilon(s) - u^*(s)) \\
&~~~ - G_1(s) \delta h(s) - G_2(s) \delta l(s)  , y(s,t) \rangle \dd s \\
&~~~ +  \langle \bar{F}(s) (u^\epsilon(s) - u^*(s)) \\
&~~~ - \bar{G}_1(s) \delta h(s) - \bar{G}_2(s) \delta l(s)  , z(s,t) \rangle \dd s \\
&~~~ + \langle \delta x(s), -Q_1(s,t) x^*(s) - \bar{Q}_1(s,t) \mathbb{E}_t[x^*(s)] \rangle \dd s \\
&~~~ + \langle \bar{H}(s) \delta x(s) + \bar{F}(s)(u^\epsilon(s) - u^*(s)) - \bar{G}_1(s) \delta h(s)  \\
&~~~ - \bar{G}_2(s) \delta l(s), y(s,t) \rangle \dd B(s)  + \langle \delta x(s),z(s,t) \rangle \dd B(s), 
\end{align*}
and
\begin{align*}
& \dd \langle \delta h(s),\phi(s,t) \rangle	 \\
& = \langle - \tilde{H}^\top(s,t) \delta h(s) - \tilde{F}^\top(s,t) \delta l(s)  \\
&~~~ - K_1^\top(s,t) (u^\epsilon(s) - u^*(s))   \\
&~~~ - K_2^\top(s,t) \mathbb{E}_t[u^\epsilon(s) - u^*(s)]   \\
&~~~ + \bar{K}_1^\top(s,t) \mathbb{E}_t[ \delta h(s)] + \bar{K}_2^\top(s,t) \mathbb{E}_t[\delta l(s)], \phi(s,t) \rangle \dd s \\
&~~~ + \langle \delta h(s), \tilde{H}(s,t) \phi(s,t) - \bar{K}_1(s,t) \mathbb{E}_t[\phi(s,t)] \\
&~~~ + G_1^\top(s) y(s,t)  + \bar{G}_1^\top(s) z(s,t)  \rangle \dd s \\
&~~~ + \langle \delta l(s), \tilde{F}(s,t) \phi(s,t) - \bar{K}_2(s,t) \mathbb{E}_t[\phi(s,t)] \\
&~~~ + G_2^\top(s) y(s,t) + \bar{G}_2^\top(s) z(s,t) \rangle \dd s \\
&~~~ + \langle \delta h(s), \tilde{F}(s,t) \phi(s,t) - \bar{K}_2(s,t) \mathbb{E}_t[\phi(s,t)] \\
&~~~ + G_2^\top(s) y(s,t) + \bar{G}_2^\top(s) z(s,t) \rangle \dd B(s)  \\
&~~~ + \langle \delta l(s), \phi(s,t) \rangle \dd B(s).
\end{align*}

In view of initial and terminal conditions of (\ref{eq_22}), we have
\begin{align*}
& \mathbb{E}_t \bigl [ \langle \delta x(T), M_1(t) x^*(T) + \bar{M}_1(t) \mathbb{E}_t[x^*(T)] \rangle 	\bigr ] \\
&  = \mathbb{E}_t  \int_t^T  \Bigl [ \langle F(s) (u^\epsilon(s) - u^*(s))  - G_1(s) \delta h(s) \\
&~~~ - G_2(s) \delta l(s)  , y(s,t) \rangle   +  \langle \bar{F}(s) (u^\epsilon(s) - u^*(s)) \\
&~~~ - \bar{G}_1(s) \delta h(s) - \bar{G}_2(s) \delta l(s)  , z(s,t) \rangle  \\
&~~~ + \langle \delta x(s), -Q_1(s,t) x^*(s) - \bar{Q}_1(s,t) \mathbb{E}_t[x^*(s)] \rangle \Bigr ] \dd s, 
\end{align*}
and
\begin{align*}	
0  
&= \mathbb{E}_t  \int_t^T \Bigl [   \langle - \tilde{H}^\top(s,t) \delta h(s) - \tilde{F}^\top(s,t) \delta l(s)  \\
& ~~~ - K_1^\top(s,t) (u^\epsilon(s) - u^*(s)) \\
&~~~  - K_2^\top(s,t) \mathbb{E}_t[u^\epsilon(s) - u^*(s)]   \\
& ~~~ + \bar{K}_1^\top(s,t) \mathbb{E}_t[ \delta h(s)] \\
&~~~ + \bar{K}_2^\top(s,t) \mathbb{E}_t[\delta l(s)], \phi(s,t) \rangle  \\
&~~~ + \langle \delta h(s), \tilde{H}(s,t) \phi(s,t)  - \bar{K}_1(s,t) \mathbb{E}_t[\phi(s,t)] \\
&~~~ + G_1^\top(s) y(s,t)  + \bar{G}_1^\top(s) z(s,t)  \rangle  \\
&~~~ + \langle \delta l(s), \tilde{F}(s,t) \phi(s,t) - \bar{K}_2(s,t) \mathbb{E}_t[\phi(s,t)] \\
& ~~~ + G_2^\top(s) y(s,t) + \bar{G}_2^\top(s) z(s,t) \rangle \Bigr ] \dd s \\
& = \mathbb{E}_t  \int_t^T \Bigl [   \langle  - K_1^\top(s,t) (u^\epsilon(s) - u^*(s)) \\
&~~~  - K_2^\top(s,t) \mathbb{E}_t[u^\epsilon(s) - u^*(s)], \phi(s,t) \rangle  \\
&~~~ + \langle \delta h(s),  G_1^\top(s) y(s,t)  + \bar{G}_1^\top(s) z(s,t)  \rangle  \\
&~~~ + \langle \delta l(s),  G_2^\top(s) y(s,t) + \bar{G}_2^\top(s) z(s,t) \rangle \Bigr ] \dd s,
\end{align*}
which implies
\begin{align}	
\label{eq_appendix_b_1}
& \mathbb{E}_t \bigl [ \langle \delta x(T), M_1(t) x^*(T) + \bar{M}_1(t) \mathbb{E}_t[x^*(T)] \rangle 	\bigr ] \\
& = \mathbb{E}_t \int_t^T  \Bigl [ \langle F(s) (u^\epsilon(s) - u^*(s))  , y(s,t) \rangle  \nonumber \\
& ~~~ +  \langle \bar{F}(s) (u^\epsilon(s) - u^*(s)) , z(s,t) \rangle  \nonumber \\
& ~~~ + \langle  - K_1^\top(s,t) (u^\epsilon(s) - u^*(s)) \nonumber \\
& ~~~ - K_2^\top(s,t) \mathbb{E}_t[u^\epsilon(s) - u^*(s)], \phi(s,t) \rangle \nonumber \\
& ~~~ + \langle \delta x(s), -Q_1(s,t) x^*(s) - \bar{Q}_1(s,t) \mathbb{E}_t[x^*(s)] \rangle \Bigr ] \dd s. \nonumber
\end{align}

On the other hand, for the objective functional of the leader in (\ref{eq_2}), we have
\begin{align*}
& J_1(t,x_0;u^\epsilon,\alpha[u^\epsilon, x_0]) - J_1(t,x_0;u^*,\alpha[u^*, x_0]) \\
& = \mathbb{E}_t \Biggl [ \int_t^T  \Bigl [ \langle x^\epsilon(s) + x^*(s), Q_1(s,t) \delta x(s) \rangle \\
&~~~ +  \langle \mathbb{E}_t [x^\epsilon(s) + x^*(s)], \bar{Q}_1(s,t) \mathbb{E}_t[\delta x(s)] \rangle  \\
&~~~ + \langle u^\epsilon(s) + u^*(s), R_1(s,t) (u^\epsilon(s) - u^*(s)) \rangle  \\
&~~~  + \langle \mathbb{E}_t[u^\epsilon(s) + u^*(s)], \bar{R}_1(s,t) \mathbb{E}_t[u^\epsilon(s) - u^*(s)] \rangle \Bigr ] \dd s\\
&~~~ + \langle x^\epsilon(T) + x^*(T), M_1(t) \delta x(T) \rangle \\
&~~~ +  \langle \mathbb{E}_t [x^\epsilon(T) + x^*(T)], \bar{M}_1(t) \mathbb{E}_t[\delta x(T)] \rangle  \Biggr ].
\end{align*}
Hence, with (\ref{eq_appendix_b_1}), we have
\begin{align}
\label{eq_appendix_b_2}
& J_1(t,x_0;u^\epsilon,\alpha[u^\epsilon, x_0]) - J_1(t,x_0;u^*,\alpha[u^*, x_0]) \\
& = \mathbb{E}_t \Biggl [ \int_t^T  \Bigl [ \langle x^\epsilon(s) + x^*(s), Q_1(s,t) \delta x(s) \rangle \nonumber \\
&~~~ +  \langle \mathbb{E}_t [x^\epsilon(s) + x^*(s)], \bar{Q}_1(s,t) \mathbb{E}_t[\delta x(s)] \rangle  \nonumber \\
&~~~ + \langle u^\epsilon(s) + u^*(s), R_1(s,t) (u^\epsilon(s) - u^*(s)) \rangle  \nonumber \\
&~~~  + \langle \mathbb{E}_t[u^\epsilon(s) + u^*(s)], \bar{R}_1(s,t) \mathbb{E}_t[u^\epsilon(s) - u^*(s)] \rangle \Bigr ] \dd s \nonumber \\
&~~~ + \langle x^\epsilon(T) + x^*(T), M_1(t) \delta x(T) \rangle \nonumber \\
&~~~ +  \langle \mathbb{E}_t [x^\epsilon(T) + x^*(T)], \bar{M}_1(t) \mathbb{E}_t[\delta x(T)] \rangle  \Biggr ] \nonumber \\
&~~~ + 2 \mathbb{E}_t \int_t^T  \Bigl [ \langle F(s) (u^\epsilon(s) - u^*(s))  , y(s,t) \rangle  \nonumber \\
&~~~ +  \langle \bar{F}(s) (u^\epsilon(s) - u^*(s)) , z(s,t) \rangle  \nonumber \\
&~~~ + \langle  - K_1^\top(s,t) (u^\epsilon(s) - u^*(s)) \nonumber \\
&~~~  - K_2^\top(s,t) \mathbb{E}_t[u^\epsilon(s) - u^*(s)], \phi(s,t) \rangle \nonumber \\
&~~~ + \langle \delta x(s), -Q_1(s,t) x^*(s) - \bar{Q}_1(s,t) \mathbb{E}_t[x^*(s)] \rangle \Bigr ] \dd s  \nonumber \\
&~~~ - 2\mathbb{E}_t \bigl [ \langle \delta x(T), M_1(t) x^*(T) + \bar{M}_1(t) \mathbb{E}_t[x^*(T)] \rangle 	\bigr ]  \nonumber \\
& = \mathbb{E}_t \Biggl [ \int_t^T  \Bigl [ \langle Q_1(s,t) (x^\epsilon(s) + x^*(s)), \delta x(s) \rangle   \nonumber \\
&~~~ + \langle  - 2 Q_1(s,t) x^*(s), \delta x(s) \rangle   \nonumber \\
&~~~ + \langle \bar{Q}_1(s,t) \mathbb{E}_t [x^\epsilon(s) + x^*(s)], \delta x(s) \rangle  \nonumber \\
&~~~ + \langle - 2 \bar{Q}_1(s,t) \mathbb{E}_t[x^*(s)], \delta x(s) \rangle  \nonumber \\
&~~~ + \langle R_1(s,t) (u^\epsilon(s) + u^*(s)), u^\epsilon(s) - u^*(s) \rangle  \nonumber \\
&~~~  + \langle \bar{R}_1(s,t)  \mathbb{E}_t[u^\epsilon(s) + u^*(s)], u^\epsilon(s) - u^*(s) \rangle   \nonumber \\
&~~~ + 2 \langle  F^\top(s) y(s,t) + \bar{F}^\top(s) z(s,t), u^\epsilon(s) - u^*(s) \rangle  \nonumber \\
&~~~ + 2 \langle -K_1(s,t) \phi(s,t) - K_2(s,t) \mathbb{E}_t[\phi(s,t)],  \nonumber \\
&~~~  u^\epsilon(s) - u^*(s) \rangle  \Bigr ] \dd s  \nonumber \\
&~~~ + \langle M_1(t)(x^\epsilon(T) + x^*(T)), \delta x(T) \rangle  \nonumber \\
&~~~ +  \langle \bar{M}_1(t) \mathbb{E}_t [x^\epsilon(T) + x^*(T)], \delta x(T) \rangle   \nonumber \\
&~~~ - 2 \langle M_1(t) x^*(T) + \bar{M}_1(t) \mathbb{E}_t[x^*(T)], \delta x(s) \rangle \Biggr ].  \nonumber
\end{align}
Then with $\widehat{Q}_1(s,t) = Q_1(s,t) + \bar{Q}_1(s,t)$ and $\widehat{M}_1(s,t) = M_1(s,t) + \bar{M}_1(s,t)$ and due to the definition of the perturbed control in (\ref{eq_6}), (\ref{eq_appendix_b_2}) can be rewritten as
\begin{align}
\label{eq_appendix_b_3}
& J_1(t,x_0;u^\epsilon,\alpha[u^\epsilon, x_0]) - J_1(t,x_0;u^*,\alpha[u^*, x_0]) \\
& = \mathbb{E}_t \Biggl [ \int_t^T \Bigl [ |\delta x(s) - \mathbb{E}_t[	\delta x(s)]|^2_{Q_1(s,t)} \nonumber \\
&~~~ + |\mathbb{E}_t[\delta x(s)]|^2_{\widehat{Q}_1(s,t)} \Bigr ] \dd s + |\mathbb{E}_t[\delta x(T)]|^2_{\widehat{M}_1(t)} \nonumber \\
&~~~ + |\delta x(T) - \mathbb{E}_t[	\delta x(T)]|^2_{M_1(t)} \nonumber \\
&~~~ + \int_t^{t+\epsilon} \Bigl [ \langle R_1(s,t) (u^\epsilon(s) + u^*(s)), u^\epsilon(s) - u^*(s) \rangle \nonumber \\
&~~~ + \langle \bar{R}_1(s,t)  \mathbb{E}_t[u^\epsilon(s) + u^*(s)], u^\epsilon(s) - u^*(s) \rangle  \nonumber \\
&~~~ + 2 \langle  F^\top(s) y(s,t) + \bar{F}^\top(s) z(s,t), u^\epsilon(s) - u^*(s) \rangle \nonumber \\
&~~~ + 2 \langle -K_1(s,t) \phi(s,t) - K_2(s,t) \mathbb{E}_t[\phi(s,t)], \nonumber \\
&~~~  u^\epsilon(s) - u^*(s) \rangle  \Bigr ] \dd s \Biggr ]. \nonumber 
\end{align}

Note that
\begin{align}
\label{eq_appendix_b_4}
& \mathbb{E}_t\Bigl [ \langle R_1(s,t) (u^\epsilon(s) + u^*(s)), u^\epsilon(s) - u^*(s) \rangle \\
&~~~ + \langle \bar{R}_1(s,t)  \mathbb{E}_t[u^\epsilon(s) + u^*(s)], u^\epsilon(s) - u^*(s) \rangle  \nonumber \\
&~~~ + 2 \langle  F^\top(s) y(s,t) + \bar{F}^\top(s) z(s,t), u^\epsilon(s) - u^*(s) \rangle \nonumber \\
&~~~ + 2 \langle -K_1(s,t) \phi(s,t) - K_2(s,t) \mathbb{E}_t[\phi(s,t)], \nonumber \\
&~~~  u^\epsilon(s) - u^*(s) \rangle  \Bigr ]	\nonumber \\
& = \mathbb{E}_t\Bigl [ \langle R_1(s,t) (u^\epsilon(s) - u^*(s)), u^\epsilon(s) - u^*(s) \rangle \nonumber \\
&~~~ + \langle \bar{R}_1(s,t)  \mathbb{E}_t[u^\epsilon(s) - u^*(s)], u^\epsilon(s) - u^*(s) \rangle  \nonumber \\
&~~~ + 2 \langle R_1(s,t) u^*(s) + \bar{R}_1(s,t) \mathbb{E}_t[u^*(s)], \nonumber \\
&~~~ u^\epsilon(s) - u^*(s) \rangle \nonumber  \\
&~~~ + 2 \langle  F^\top(s) y(s,t) + \bar{F}^\top(s) z(s,t), u^\epsilon(s) - u^*(s) \rangle \nonumber  \\
&~~~ + 2 \langle -K_1(s,t) \phi(s,t) - K_2(s,t) \mathbb{E}_t[\phi(s,t)], \nonumber  \\
&~~~  u^\epsilon(s) - u^*(s) \rangle  \Bigr ] \nonumber \\
& = \mathbb{E}_t \Bigl [|(u^\epsilon(s) - u^*(s)) - \mathbb{E}_t[u^\epsilon(s) - u^*(s)] |^2_{R_1(s,t)}  \nonumber \\
&~~~ + |\mathbb{E}_t[u^\epsilon(s) - u^*(s)] |^2_{\widehat{R}_1(s,t)} \nonumber \\
&~~~ + 2 \langle \lambda(s,t),u^\epsilon(s) - u^*(s) \rangle \Bigr ], \nonumber 
\end{align}
where
\begin{align*}
\lambda(s,t) &= 	R_1(s,t) u^*(s) + \bar{R}_1(s,t) \mathbb{E}_t[u^*(s)] \\
&~~~ + F^\top(s) y(s,t) + \bar{F}^\top(s) z(s,t) \\
&~~~ -K_1(s,t) \phi(s,t) - K_2(s,t) \mathbb{E}_t[\phi(s,t)].
\end{align*}
Due to the definition of $u^*$ in (\ref{eq_24}), boundedness of the coefficients, and the continuity of the FBSDEs in (\ref{eq_20}) and (\ref{eq_22}), we have
\begin{align}
\label{eq_appendix_b_5}
& \liminf_{s \downarrow t} \mathbb{E}[|\lambda(s,t)|^2] = 0.
\end{align}
Hence, (\ref{eq_appendix_b_4}) can be rewritten as follows:
\begin{align*}
& \mathbb{E}_t \Bigl [|(u^\epsilon(s) - u^*(s)) - \mathbb{E}_t[u^\epsilon(s) - u^*(s)] |^2_{R_1(s,t)}  \\
&~~~ + |\mathbb{E}_t[u^\epsilon(s) - u^*(s)] |^2_{\widehat{R}_1(s,t)} \\
&~~~ + 2 \langle \lambda(s,t),u^\epsilon(s) - u^*(s) - \mathbb{E}_t[u^\epsilon(s) - u^*(s)] \\
&~~~ + \mathbb{E}_t[u^\epsilon(s) - u^*(s)] \rangle \Bigr ]	 \\
& =  \mathbb{E}_t \Biggl [ \Bigl |R_1^{\frac{1}{2}}(s,t)(u^\epsilon(s) - u^*(s)) \\
&~~~ - R_1^{\frac{1}{2}}(s,t)\mathbb{E}_t[u^\epsilon(s) - u^*(s)]  + R_1^{-\frac{1}{2}}(s,t) \lambda(s,t) \Bigr |^2  \\
&~~~ + \Bigl |\widehat{R}_1^{\frac{1}{2}}(s,t)\mathbb{E}_t[u^\epsilon(s) - u^*(s)] + \widehat{R}_1^{-\frac{1}{2}}(s,t) \mathbb{E}_t[\lambda(s,t)] \Bigr |^2 \\
&~~~ - \bigl | \lambda(s,t) \bigr |^2_{R_1^{-1}(s,t)} - \bigl | \mathbb{E}_t[\lambda(s,t)] \bigr |^2_{\widehat{R}_1^{-1}(s,t)} \Biggr ].
\end{align*}

This and (\ref{eq_appendix_b_3}) lead to
\begin{align*}
& J_1(t,x_0;u^\epsilon,\alpha[u^\epsilon, x_0]) - J_1(t,x_0;u^*,\alpha[u^*, x_0]) \\
& = \mathbb{E}_t \Biggl [ \int_t^T \Bigl [ |\delta x(s) - \mathbb{E}_t[	\delta x(s)]|^2_{Q_1(s,t)} \\
&~~~ + |\mathbb{E}_t[\delta x(s)]|^2_{\widehat{Q}_1(s,t)} \Bigr ] \dd s + |\mathbb{E}_t[\delta x(T)]|^2_{\widehat{M}_1(t)} \\
&~~~ + |\delta x(T) - \mathbb{E}_t[	\delta x(T)]|^2_{M_1(t)} \\
&~~~ + \int_t^{t+\epsilon}	\Bigl [  \Bigl |R_1^{\frac{1}{2}}(s,t)(u^\epsilon(s) - u^*(s)) \\
&~~~ - R_1^{\frac{1}{2}}(s,t)\mathbb{E}_t[u^\epsilon(s) - u^*(s)]  + R_1^{-\frac{1}{2}}(s,t) \lambda(s,t) \Bigr |^2  \\
&~~~ + \Bigl |\widehat{R}_1^{\frac{1}{2}}(s,t)\mathbb{E}_t[u^\epsilon(s) - u^*(s)] + \widehat{R}_1^{-\frac{1}{2}}(s,t) \mathbb{E}_t[\lambda(s,t)] \Bigr |^2 \\
&~~~ - \bigl | \lambda(s,t) \bigr |^2_{R_1^{-1}(s,t)} - \bigl | \mathbb{E}_t[\lambda(s,t)] \bigr |^2_{\widehat{R}^{-1}_1(s,t)} \Bigr ] \dd s \Biggr ] \\
& \geq - \mathbb{E}_t  \int_t^{t+\epsilon} \Bigl [ \bigl | \lambda(s,t) \bigr |^2_{R_1^{-1}(s,t)} + \bigl | \mathbb{E}_t[\lambda(s,t)] \bigr |^2_{\widehat{R}_1^{-1}(s,t)}  \Bigr ] \dd s ,
\end{align*}
which, together with (\ref{eq_appendix_b_5}), implies
\begin{align*}
& \liminf_{\epsilon \downarrow 0} \frac{1}{\epsilon} \Bigl (J_1(t,x_0;u^\epsilon,\alpha[u^\epsilon, x_0]) \\
&~~~ - J_1(t,x_0;u^*,\alpha[u^*, x_0]) \Bigr ) \geq 0
\end{align*}
Hence, in view of Definition \ref{Definition_2}(ii), we have the desired result. This completes the proof. $\hfill{\square}$

\section{Explicit Expressions of the Nonsymmetric Coupled RDEs in (\ref{eq_34})}\label{Appendix_C}
In this appendix, we provide the explicit expressions of nonsymmetric coupled RDEs in (\ref{eq_34}). 

In (\ref{eq_34}), we have
\begin{align*}
\begin{cases}
- \frac{\dd \mathcal{P}(s,t)}{\dd s} = \Lambda_1(s,t, \mathcal{P}(s,t), \mathcal{Z}(s,t)),~ s \in [t,T]\\
- \frac{\dd \mathcal{Z}(s,t)}{\dd s} = \Lambda_2(s,t, \mathcal{P}(s,t), \mathcal{Z}(s,t)),~ s \in [t,T]\\
\mathcal{P}(T,t) = \mathcal{M}(t),~ \mathcal{Z}(T,t) = \bar{\mathcal{M}}(t)		\\
	\det (I - \mathcal{P}(s,t) \bar{\mathcal{D}}(s))  \neq 0 \\
	\det \bigl (\widehat{\mathcal{R}}(t,t)  + \bar{\mathcal{B}}^\top(t) \bigl [ I - \mathcal{P}(t,t) \bar{\mathcal{D}}(t) \bigr ]^{-1} \mathcal{P}(t,t)  \bar{\mathcal{B}}(t) \bigr) \neq 0\\
		\widehat{\mathcal{P}}(s,t) = \mathcal{P}(s,t) + \mathcal{Z}(s,t),
\end{cases}	
\end{align*}
where
\begin{align}
\label{eq_appendix_c_1}
& \Lambda_1(s,t,\mathcal{P}(s,t),\mathcal{Z}(s,t))	\\
& = \mathcal{A}_1^\top (s) \mathcal{P}(s,t) + \mathcal{P}(s,t) \mathcal{A}_1(s)  + \mathcal{Q}(s,t)  \nonumber \\
&~~~ + \mathcal{P}(s,t) \mathcal{C}(s) \mathcal{P}(s,t) + \Bigl [ \bar{\mathcal{A}}_1^\top(s,t) + \mathcal{P}(s,t) \mathcal{D}(s) \Bigr ]  \nonumber \\
&~~~ \times \Bigl [ I - \mathcal{P}(s,t) \bar{\mathcal{D}}(s) \Bigr ]^{-1} \mathcal{P}(s,t) \Bigl ( \bar{\mathcal{A}}_1(s,t)  \nonumber \\
&~~~ - \bar{\mathcal{B}}(s) \Bigl [ \widehat{\mathcal{R}}(s,s) + \bar{\mathcal{B}}^\top(s) \bigl [ I - \mathcal{P}(s,s) \bar{\mathcal{D}}(s) \bigr ]^{-1} \mathcal{P}(s,s)  \bar{\mathcal{B}}(s) \Bigr ]^{-1}
\nonumber \\
&~~~ \times \Bigl [ \mathcal{B}^\top(s) \widehat{\mathcal{P}}(s,s) +\bar{\mathcal{B}}^\top(s) \bigl [ I - \mathcal{P}(s,s) \bar{\mathcal{D}}(s) \bigr ]^{-1} \nonumber \\
&~~~ \times  \mathcal{P}(s,s) (\widehat{\mathcal{A}}_2(s,s) + \bar{\mathcal{C}}(s) \widehat{\mathcal{P}}(s,s) )\Bigr ]
 + \bar{\mathcal{C}}(s) \mathcal{P}(s,t) \Bigr) \nonumber \\
&~~~ - \Bigl [ \mathcal{G}^\top(s,t) + \mathcal{P}(s,t) \mathcal{B}(s)  \Bigr ] \nonumber \\
&~~~ \times \Bigl [ \widehat{\mathcal{R}}(s,s) + \bar{\mathcal{B}}^\top(s) \bigl [ I - \mathcal{P}(s,s) \bar{\mathcal{D}}(s) \bigr ]^{-1}  \mathcal{P}(s,s)  \bar{\mathcal{B}}(s) \Bigr ]^{-1}  
\nonumber \\
&~~~ \times \Bigl [ \mathcal{B}^\top(s) \widehat{\mathcal{P}}(s,s)  +\bar{\mathcal{B}}^\top(s) \bigl [ I - \mathcal{P}(s,s) \bar{\mathcal{D}}(s) \bigr ]^{-1}   \nonumber \\
&~~~ \times \mathcal{P}(s,s) (\widehat{\mathcal{A}}_2(s,s) + \bar{\mathcal{C}}(s) \widehat{\mathcal{P}}(s,s) )\Bigr ], \nonumber
\end{align}
and
\begin{align}
\label{eq_appendix_c_2}		
& \Lambda_2(s,t,\mathcal{P}(s,t),\mathcal{Z}(s,t))	\\
& = \mathcal{A}_1^\top(s,t) \mathcal{Z}(s,t) + \mathcal{A}_2^\top(s,t) \widehat{\mathcal{P}}(s,t) + \mathcal{P}(s,t) \mathcal{A}_2(s,t) \nonumber \\
&~~~ + \mathcal{Z}(s,t) \widehat{\mathcal{A}}_1(s,t) + \bar{\mathcal{Q}}(s,t) \nonumber \\
&~~~  + \mathcal{P}(s,t) \mathcal{C}(s) \mathcal{Z}(s,t)   + \mathcal{Z}(s,t) \mathcal{C}(s) \widehat{\mathcal{P}}(s,t) \nonumber \\
&~~~ + \Bigl [ \bar{\mathcal{A}}_1^\top(s,t) + \mathcal{P}(s,t) \mathcal{D}(s) \Bigr ] \Bigl [ I - \mathcal{P}(s,t) \bar{\mathcal{D}}(s) \Bigr ]^{-1}  	\mathcal{P}(s,t) \nonumber \\
&~~~ \times \Bigl [\bar{\mathcal{A}}_2(s,t) + \bar{\mathcal{C}}(s) \mathcal{Z}(s,t) \Bigr ] \nonumber \\
&~~~ + \Bigl [ \bar{\mathcal{A}}_2^\top(s,t) + \mathcal{Z}(s,t) \mathcal{D}(s) \Bigr ] \Bigl [ I - \mathcal{P}(s,t) \bar{\mathcal{D}}(s) \Bigr ]^{-1} \mathcal{P}(s,t) \nonumber \\
&~~~ \times \Bigl ( \widehat{\mathcal{A}}_2(s,t) + \bar{\mathcal{C}}(s) \widehat{\mathcal{P}}(s,t)  - \bar{\mathcal{B}}(s)  \nonumber \\
&~~~  \times \Bigl [ \widehat{\mathcal{R}}(s,s)  + \bar{\mathcal{B}}^\top(s) \bigl [ I - \mathcal{P}(s,s) \bar{\mathcal{D}}(s) \bigr ]^{-1}  \mathcal{P}(s,s)  \bar{\mathcal{B}}(s) \Bigr ]^{-1} \nonumber \\
&~~~ \times \Bigl [ \mathcal{B}^\top(s) \widehat{\mathcal{P}}(s,s) + \bar{\mathcal{B}}^\top(s) \bigl [ I - \mathcal{P}(s,s) \bar{\mathcal{D}}(s) \bigr ]^{-1}   \nonumber \\
&~~~ \times \mathcal{P}(s,s) (\widehat{\mathcal{A}}_2(s,s) + \bar{\mathcal{C}}(s) \widehat{\mathcal{P}}(s,s) )\Bigr ] 
   \Bigr ) \nonumber \\
&~~~  - \Bigl [ \bar{\mathcal{G}}^\top(s,t) + \mathcal{Z}(s,t) \mathcal{B}(s) \Bigr ] \nonumber \\
&~~~ \times \Bigl [ \widehat{\mathcal{R}}(s,s) + \bar{\mathcal{B}}^\top(s) \bigl [ I - \mathcal{P}(s,s) \bar{\mathcal{D}}(s) \bigr ]^{-1}  \mathcal{P}(s,s)  \bar{\mathcal{B}}(s) \Bigr ]^{-1}
\nonumber \\
&~~~ \times \Bigl ( \mathcal{B}^\top(s) \widehat{\mathcal{P}}(s,s) +  \bar{\mathcal{B}}^\top(s) \bigl [ I - \mathcal{P}(s,s) \bar{\mathcal{D}}(s) \bigr ]^{-1}  \nonumber \\
&~~~ \times  \mathcal{P}(s,s) (\widehat{\mathcal{A}}_2(s,s) + \bar{\mathcal{C}}(s) \widehat{\mathcal{P}}(s,s) )\Bigr ). \nonumber
\end{align}
Hence, we have
\begin{align}
\label{eq_appendix_c_3}
\begin{cases}
	\widehat{\mathcal{P}}(s,t) = \mathcal{P}(s,t) + \mathcal{Z}(s,t) \\
	\frac{\dd \widehat{\mathcal{P}}(s,t)}{\dd s} = 
	\Lambda_3(s,t,\mathcal{P}(s,t),\mathcal{Z}(s,t)),~ s \in [t,T]
	\\
	\widehat{\mathcal{P}}(T,t) = \mathcal{M}(t) + \bar{\mathcal{M}}(t),
\end{cases}	
\end{align}
where the explicit expression of $\Lambda_3$ can be obtained from (\ref{eq_appendix_c_1}) and (\ref{eq_appendix_c_2}).

\bibliographystyle{IEEEtran}       
\bibliography{researches_1.bib}

\end{document}